\documentclass[a4paper,french,12pt]{article}

\usepackage[utf8]{inputenc}
\usepackage[english]{babel}  
\usepackage[T1]{fontenc}      
\usepackage{geometry}         
\usepackage{amsfonts}
\usepackage{amsmath}
\usepackage{amssymb}
\usepackage[usenames, dvipsnames]{xcolor}
\usepackage[pdftex]{graphicx}
\usepackage{graphics}
\usepackage{float}
\usepackage[all]{xy}
\usepackage{color}
\usepackage{fourier}
\usepackage{tikz}
\usepackage{amsthm}
\usetikzlibrary{patterns}

\def\R{\mathop{\mathbb{R}}\nolimits}
\def\Z{\mathop{\mathbb{Z}}\nolimits}
\def\B{\mathop{\mathbb{B}}\nolimits}

\def\ksi {\mathop{\xi}\nolimits}

\def\* {\mathop{\otimes}\nolimits}
\def\+{\mathop{\oplus}\nolimits}

\newtheorem{theorem}{Theorem}[section]
\newtheorem{corollary}[theorem]{Corollary}
\newtheorem{lemma}[theorem]{Lemma}
\newtheorem{proposition}[theorem]{Proposition}
\newtheorem{definition}{Definition}[section]
\newtheorem*{theorem*}{Theorem}
\newtheorem*{lemma*}{Lemma}
\newtheorem*{proposition*}{Proposition}
\newtheorem*{corollary*}{Corollary}
\theoremstyle{remark}
\newtheorem*{remark}{Remark}

\def\Dast{D_{\ast}\left(1/3\right)}
\def\Dtw{D\left(1/2\right)}
\def\Dth{D\left(1/3\right)}
\newcommand{\Alpha}{0.4}
\newcommand{\Aalpha}{0.60}

\def\les{\lesssim}

\def\eps{\varepsilon}
\def\Dastp#1{D^{#1}_{\ast}\left(1/3\right)}

\begin {document}
\title{New bounds for the inhomogenous Burgers and  the Kuramoto-Sivashinsky equations}
\author{Michael Goldman\footnote{LJLL, Universit\'e Paris Diderot,  CNRS, UMR 7598,  France, email: goldman@math.univ-paris-diderot.fr}, Marc Josien\footnote{Ecole Polytechnique, Palaiseau, France, email: marc.josien@polytechnique.edu}  and Felix Otto\footnote{Max Planck Institute for Mathematics in the Sciences, Leipzig, Germany, email: otto@mis.mpg.de}}

\maketitle

\begin{abstract}
We give a substantially simplified proof of the near-optimal
estimate on the Kuramoto-Sivashinsky equation from \cite{Otto}, at the same time slightly
improving the result. The result in \cite{Otto} relied on two ingredients: a regularity
estimate for capillary Burgers and an a novel priori estimate for the {\it inhomogeneous inviscid}
Burgers equation, which works out that in many ways the {\it conservative transport} nonlinearity
acts as a {\it coercive} term. It is the proof of the second ingredient that we substantially
simplify by proving a modified K\'arm\'an-Howarth-Monin identity for solutions of the  inhomogeneous inviscid Burgers equation. We show that this provides a new interpretation of the results obtained in \cite{Golse}.
\end{abstract}

\section{Introduction}
	\subsection{The Kuramoto-Sivashinsky equation}
\begin{figure}[ht]\label{chaos}
			\begin{center}
			\includegraphics[width=16cm]{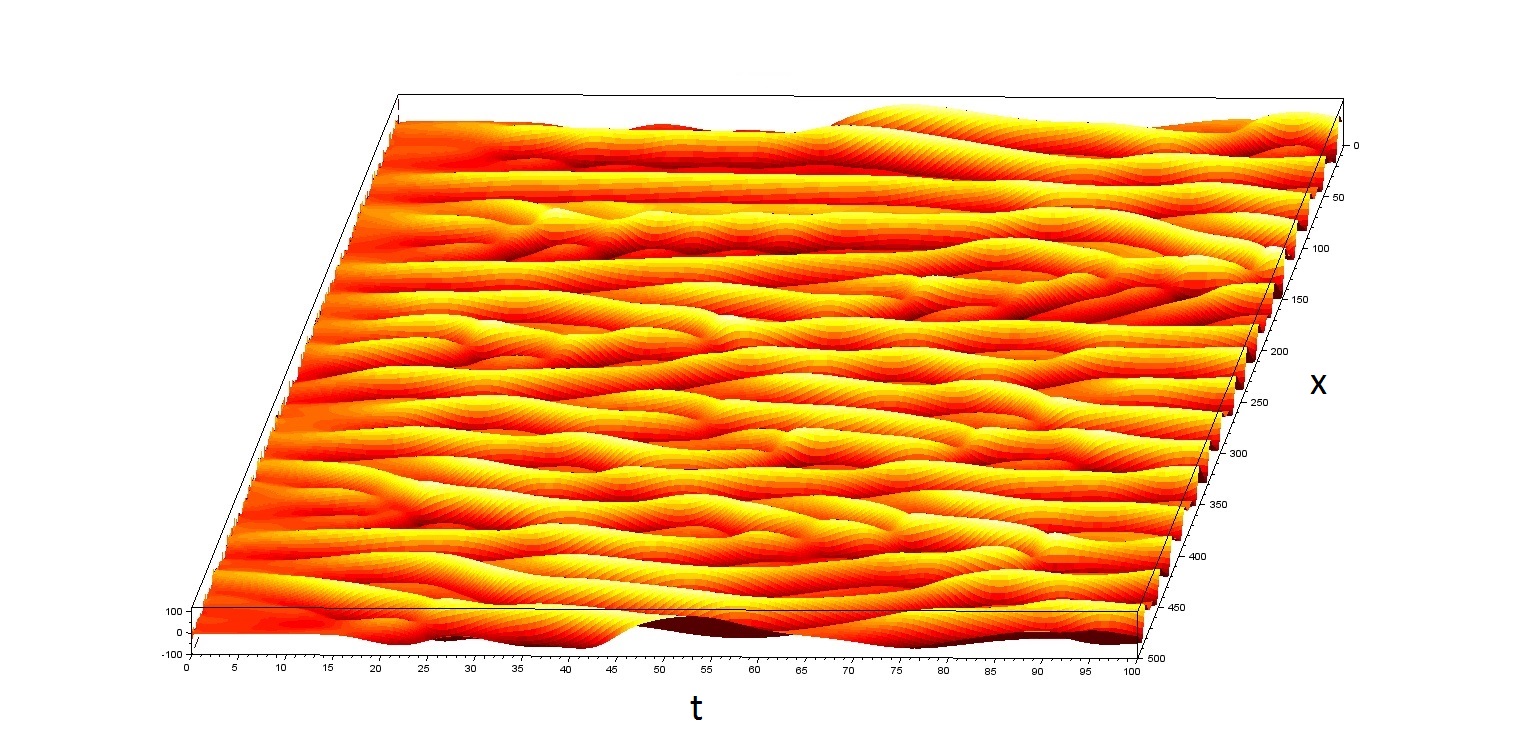} \caption{Chaotic behavior of $u$}
			\end{center}
		\end{figure}

		We consider the one-dimensional Kuramoto-Sivashinsky equation:
		\begin{equation}\label{Kuramoto} \tag{K-S}
			\partial_t u +  u \partial_x u + \partial_x^2u+\partial_x^4u=0.
		\end{equation}
		This equation appears in many physical contexts, in particular in the  modeling of surface evolutions. Sivashinsky used it to describe flame fronts \cite{Sivashinsky}, wavy flow of viscous liquids on inclined planes \cite{Syva2} and crystal growth \cite{Golovin}.
		Although  the solutions of \eqref{Kuramoto} are smooth and even analytic \cite{Jolly},  they display a chaotic behavior for sufficiently large systems size $L$ (see \cite{Hyman} and Figure \ref{chaos}). 
		The structure of the Kuramoto-Sivashinsky equation has some similarities with the Navier-Stokes equation. Therefore, it is sometimes possible to apply similar techniques
		to study both equations (see \cite{Ramos}).\\

		For a given system size $L>0$, we will consider $L$-periodic solutions of \eqref{Kuramoto}. Since the spatial average $\int_0^L u(t,x)dx$ is constant over time, and since the equation is invariant under the Galilean transformation:
		\begin{align*}
			t=t', &&x=x'+Ut,&&u=u'+U,
		\end{align*}
		it is not restrictive to assume that $\int_0^Lu(t,x)dx=0$ for all $t\geq0$.\\
		We can artificially cut the equation in two parts and consider separately the two mechanisms involved in \eqref{Kuramoto}:
		\begin{align}
			\partial_tu+\partial_x^2u+\partial_x^4u=0, \label{phénomène1}\\
			\partial_tu+u\partial_xu=0. \label{phénomène2}
		\end{align}
		The first equation \eqref{phénomène1} is linear and  can be seen in Fourier space as:
		\begin{equation}\label{linear part}
			\partial_t\mathcal{F}u=(\ksi^2-\ksi^4)\mathcal{F}u.
		\end{equation}
		The fourth partial derivative term $\partial_x^4u$ decreases the short wavelength part of the energy spectrum whereas the second derivative term $\partial_x^2u$ amplifies the long wavelength part.
		The second equation \eqref{phénomène2} corresponds to  Burgers equation. It is nonlinear and develops shocks in finite time for non-trivial initial data.
		Nevertheless, as we will see later, this term has some mild regularizing effect. It is worth mentioning that for \eqref{phénomène2}, the energy
		\begin{equation*}
			\int_0^L u^2 dx,
		\end{equation*}
		is conserved. Therefore, one can intuitively say that in \eqref{Kuramoto} the linear terms transport the energy from long wavelengths to short ones. Numerical simulations suggest (see the article of Wittenberg and Holmes \cite{Holmes}) that the time-averaged power spectrum
		\begin{equation*}
			\lim_{T \rightarrow \infty} (LT)^{-1} \int_0^T \left| \mathcal{F}(u)(t,\ksi) \right|^2dt
		\end{equation*}
		is independent of $L$ for $L\gg1$. Moreover, this quantity is independent of   $\left|\ksi\right|$ and $L\gg 1$ in the long wavelength regime $L^{-1}\ll |\ksi|\ll 1$ and decays exponentially  in the short wavelength regime $\left|\ksi\right| \gg 1$. In line with this, numerical simulations suggest that for all $\alpha \geq 0$:
		\begin{equation*}
			\limsup_{T\rightarrow \infty} \, (LT)^{-1} \int_0^T\int_0^L \left(\left|\partial_x \right|^\alpha u\right)^2dxdt=O(1).
		\end{equation*}
		This conjecture is supported by a universal bound on all stationary periodic solutions of \eqref{Kuramoto} with mean $0$, due to Michelson \cite{Michelson}.

	\subsection{Known bounds}
		A first energy bound was obtained in the 80's by Nicolaenko, Scheurer and Temam \cite{NST}, who established by the ``background flow method'' that for every odd (in space) solution $u$ of \eqref{Kuramoto}:
		\begin{equation*}
			\limsup_{t\rightarrow \infty} \left(\frac{1}{L}\int_0^L u^2dx\right)^{1/2} = O\left(L^p\right),
		\end{equation*}
		with $p=2$. This has been later  generalized by Goodman \cite{Goodman}, and Bronski and Gambill \cite{Bronski} and  improved to $p=1$. Using an entropy method, Giacomelli and the third author \cite{Giacomelli} improved this result by showing that:
		\begin{equation*}
			\limsup_{t\rightarrow \infty} \left(\frac{1}{L}\int_0^L u^2dx\right)^{1/2} = o(L).
		\end{equation*}
		The proof is based on the fact that the dispersion relation $\ksi^2-\ksi^4$ in \eqref{linear part} vanishes for $\ksi \rightarrow 0$ and it implies  that for every $\alpha \in [0,2]$, we have:
		\begin{equation*}
			\limsup_{T\rightarrow \infty} \left(\frac{1}{TL}\int_0^T \int_0^L \left(\left|\partial_x\right|^\alpha u\right)^2dxdt\right)^{1/2} = o(L),
		\end{equation*}
		by using the energy identity,
		\begin{equation*}
			\partial_t \int_0^L\left(u(t,x)\right)^2dx=\int_0^L \left(\partial_x u(t,x)\right)^2dx-\int_0^L \left( \partial_x^2 u(t,x)\right)dx.
		\end{equation*}
		  In a more recent paper \cite{Otto}, the third author proved that, for all $\alpha \in (1/3,2)$,
		  \begin{equation}\label{firstOtto}
		  	\limsup_{T\rightarrow \infty} \left(\frac{1}{TL}\int_0^T \int_0^L \left(\left|\partial_x\right|^\alpha u\right)^2dxdt\right)^{1/2} = O\left(\ln^{5/3+}(L)\right) ,\footnote{we use the notation $O\left(\ln^{5/3+}(L)\right)$ to indicate that for every $\kappa>5/3$, there exists $c>0$ such that the left-hand side is bounded by $ c \ln^\kappa(L)$.}
		  \end{equation}
		  by using two ingredients: an {\it a priori} estimate for the capillary Burgers equation  
			$\partial_t u +u\partial_x u +\partial_x^4u=|\partial_x|g$ and an {\it a priori} estimate for  the inhomogeneous  Burgers equation, that is $\partial_t u +u\partial_x u =|\partial_x|g$. More precisely, the result of \cite{Otto} states that, for every solution $u$ of \eqref{Kuramoto},
		  \begin{equation}\label{rsltOtto}
		  \left\|u\right\|_{\B^{1/3}_{3,3}} = O\left(\ln^{5/3+}(L) \right),
		  \end{equation}
		  where  $\left\| \cdot \right\|_{\B^{s}_{p,r}}$ denotes a Besov norm (see the appendix).
		  
		  \subsection{Main result}
		  	In this paper, we improve and simplify the result of the third author by showing that:
		 \begin{theorem}
		Let $L>2$. For $u$ a smooth $L$-periodic solution with zero average of the equation
				\begin{equation*}
					\partial_tu+u\partial_xu+\partial_x^2u+\partial_x^4u=0,				
				\end{equation*}	
				there holds 
		\begin{equation}\label{MainIntro}
		  	\limsup_{T\rightarrow \infty} \left( \sup_{h>0}\frac{1}{LT}\int_0^T \int_0^L \frac{\left|u(t,x+h)-u(t,x)\right|^3}{h}dxdt \right)^{1/3}= O\left(\ln^{1/2+}(L)\right).
		  \end{equation}
		\end{theorem}
			This result is indeed slightly stronger than the previous one, since by \eqref{logarithm}, it implies an improvement of the exponent in \eqref{rsltOtto} from $5/3+$ to  $5/6+$. However, this is not the main contribution of this paper. It 
		  is rather a simplified proof of the a priori estimate for inhomogeneous  Burgers equation, which was one of the main tool for proving \eqref{rsltOtto}. 
		   For this purpose, we derive a modified K\'arm\'an-Howarth-Monin formula (see \eqref{KHMstrat}). We also show how the proof of  Golse and Perthame \cite{Golse} 
		   (based on the kinetic formulation of Burgers equation) of a similar estimate for the homogeneous Burgers equation  can be reinterpreted in this light. 
		   Since we work with slightly different Besov norms compared to \cite{Otto}, 
		   we need also to adapt most of the other steps to get \eqref{MainIntro}. Besides Proposition \ref{thm2}, which we borrow directly from \cite{Otto},
		   we give here self-contained proofs. \\
		  
		 The structure of the paper is the following:  In Section \ref{sec:Main}, we enunciate the main theorem and give the structure of the proof. It has several ingredients: a Besov estimate for the inhomogeneous inviscid  equation (Proposition \ref{thm1}), a regularity estimate for the capillary Burgers equation (Proposition \ref{thm2}) and an 
		 inverse estimate for Besov norms on solutions of \eqref{Kuramoto} (Proposition \ref{thm3}). The following sections (\textit{i.e.} Section \ref{sec:proofmain}, \ref{sec:proofburgers} and \ref{sec:increaser}) are devoted to the proofs. In the appendix, we recall definitions and a few classical results regarding Besov spaces.

\section*{General notations}
	
		We denote by $D^h $ the finite-difference operator $D^h :u \mapsto u(x+h)-u(x)$, by  $L^p$ the space $L^p\left([0,T] \times [0,L]\right)$ and for $k\in \mathbb{N}, L>0$, by  $\mathcal{C}_L^k= \left\{f \in \mathcal{C}^k(\R), f \text{ is } L-\text{periodic} \right\}$.\\	
		For an  $L$-periodic function $u$, the spatial Fourier transform is defined by:
		\begin{equation*}
			\mathcal{F}(u)(\ksi)=L^{-1}\int_0^L\exp(-i\ksi x)u(x)dx
		\end{equation*}
		and for a Schwartz function $\phi$:
		\begin{equation*}
		\mathcal{F}(\phi)(\ksi)=\int_{\R} \exp(-i\ksi x)\phi(x)dx.
		\end{equation*}
		For  $v\in \R$, we let $v_+=\max(v,0)$ (and similarly, $v_-= \max(-v,0)$).
\section{Main theorem and structure of the proof}\label{sec:Main}
	In this section, we state the main theorem and the results on which it is based (see the appendix for the definition and main properties of Besov spaces).
	\subsection{Main theorem}
			\begin{theorem} \label{maintheorem} 
			 Let $L>2$. For a smooth $L$-periodic solution  $u$ with zero average of \eqref{Kuramoto},	there holds 
				\begin{equation}\label{Main}
					\left\|u \right\|_{\B^{1/3}_{3,\infty}}+\left\|u \right\|_{\B^2_{2,2}} =O(  \ln^{1/2+}(L)).
				\end{equation}				
		\end{theorem}
		From this theorem, we derive  by interpolation \eqref{Interpolation} the following corollary:
		\begin{corollary}\label{corol}
			 Let $L>2$. For a smooth $L$-periodic solution $u$ with zero average  of \eqref{Kuramoto}
				and for indices $p$, $s$ and $r$ related by
				\begin{align*}
				 p \in [1/3,2], && s=10/p-3, && 1/r = 3(1/p-1/3),
				\end{align*}
				we have
				\begin{equation*}
					\left\|u \right\|_{\B^{s}_{p,r}} =O\left(  \ln^{1/2+}(L)\right).
				\end{equation*}
		\end{corollary}
		
	\subsection{Structure of the proof}
\begin{figure}[H]\caption{Strategy of the proof} \label{bellefigure}
\begin{tikzpicture}[scale=13]
	\coordinate (O) at (0,0) ;
	\coordinate (A1) at (3/10,0) ;
	\coordinate (B1) at (3/10,1);
	\coordinate (A2) at (1/3,0);
	\coordinate (A21)  at (1/3,-0.04);
	\coordinate (B2) at (1/3,1);
	\coordinate (A3) at (2/5,0);
	\coordinate (B3) at (2/5,1);
	\coordinate (A4) at (1/2,0);
	\coordinate (C4) at (0,1/2);
	
	\coordinate (A5) at (3/5,0);
	\coordinate (B5) at (3/5,1);
	\coordinate (A6) at (2/3,0);
	\coordinate (A61) at (2/3,-0.04);
	\coordinate (B6) at (2/3,1);
	\coordinate (A7) at (4/5,0);
	\coordinate (B7) at (4/5,1);
	
	\coordinate (delta) at (2/5,2/5);
	\coordinate (Beta) at (4/5,4/5);
	\coordinate (J) at (1/2,1/2);

	\newcommand{\point}[2]%
	{\draw (#1) node {$\bullet$};
	\draw (#1) node [below] {$#2$};
	}
	
	\newcommand{\droite}[3]%
	{\draw[line width=0.5mm,#1] (#2)--(#3);
	}
	
	\draw (A21) node {$1/3$};
	\draw (A2) node {$\bullet$};
	\point {A3}{2/5};
	\point {A4}{1/2};
	\point {A5}{3/5};
	\point{A7}{4/5};
	
	\draw [thick,->,>=stealth] (0,0)--(0,1);
	\draw [thick,->,>=stealth] (0,0)--(1,0);

	\draw [dashed] (A4)--(J);
	\draw [dashed] (C4)--(J);
	\draw (C4) node [left] {$1/2$};
	\draw (J) node {$\bullet$};

	\draw [dashed, line width=0.5mm, RoyalBlue] (0,0)--(1/3,1/3);
	\draw [line width=0.5mm,RoyalBlue] (1/3,1/3)--(4/5,4/5);
	\draw [dashed, line width=0.5mm, RoyalBlue] (4/5,4/5)--(1,1);
	
	\draw [<->,line width=1mm, RoyalBlue] (1/3,1/3)--(1/2,1/2);
	\draw [<->, line width=1mm, RoyalBlue] (1/2,1/2)--(4/5,4/5);

	\draw (\Alpha,\Alpha) node {$\bullet$};
	
	\draw (\Aalpha,\Aalpha) node {$\bullet$};

	\draw (Beta) node {$\bullet$};

	\draw [line width=1 mm, OliveGreen,->] (1/3,1/3)--(1/3,0);
	\draw (1/3,1/3) node {$\bullet$};
	
	\draw [line width=1mm,<-] (delta) ..controls +(0,0.3) and +(0,0.3).. (Beta);
	
	\draw [line width=1mm, ->,Red] (1/3,0) ..controls +(0.2,0) and +(0,-0.3).. (\Aalpha,\Aalpha);
	\draw [line width=1mm, ->,Red] (1/3,0) ..controls +(0.1,0) and +(0,-0.3).. (\Alpha,\Alpha);
	
	\draw [line width = 1mm,->, Purple] (1/2,1/2) ..controls +(-0.3,0.3) and (-0.1,0.3).. (1/3,0);

\coordinate (E1) at (0,-0.07);
\coordinate (E2) at (0.75,-0.07);
\coordinate (E3) at (0.75,-0.29);
\coordinate (E4) at (0,-0.29);
\draw (E1)--(E2)--(E3)--(E4)--cycle;

\draw (E1) node [below right] {\textbf{Legend:}};

\draw (0,-0.14) node [right,RoyalBlue] {$\bullet$};
\draw (0.05,-0.14) node [right] {Interpolation};

\draw (0,-0.17) node [right, Purple] {$\bullet$};
\draw (0.05,-0.17) node [right] {Energy estimate};

\draw (0,-0.20) node [right, Black] {$\bullet$};
\draw (0.05,-0.20) node [right] {Higher regularity};

\draw (0,-0.23) node [right, OliveGreen] {$\bullet$};
\draw (0.05,-0.23) node [right] {Increase of the index $r$};

\draw (0,-0.26) node [right, Red] {$\bullet$};
\draw (0.05,-0.26) node [right] {Estimate of the inhomogeneous Burgers equation};

	\draw (1/3,1/3) node {$\bullet$};
	\draw (Beta) node {$\bullet$};
	\draw (delta) node {$\bullet$};
	\draw (\Aalpha,\Aalpha) node {$\bullet$};
	\draw (\Alpha,\Alpha) node {$\bullet$};
	\draw (J) node {$\bullet$};
	\draw (A21) node {$1/3$};
	\draw (A2) node {$\bullet$};
	\point {A3}{2/5};
	\point {A4}{1/2};
	\point {A5}{3/5};
	\point{A7}{4/5};
	
	\draw (0,1) node [left] {$1/r$};	
	\draw (1,0) node [below] {$1/p$};															
\end{tikzpicture}
\end{figure}
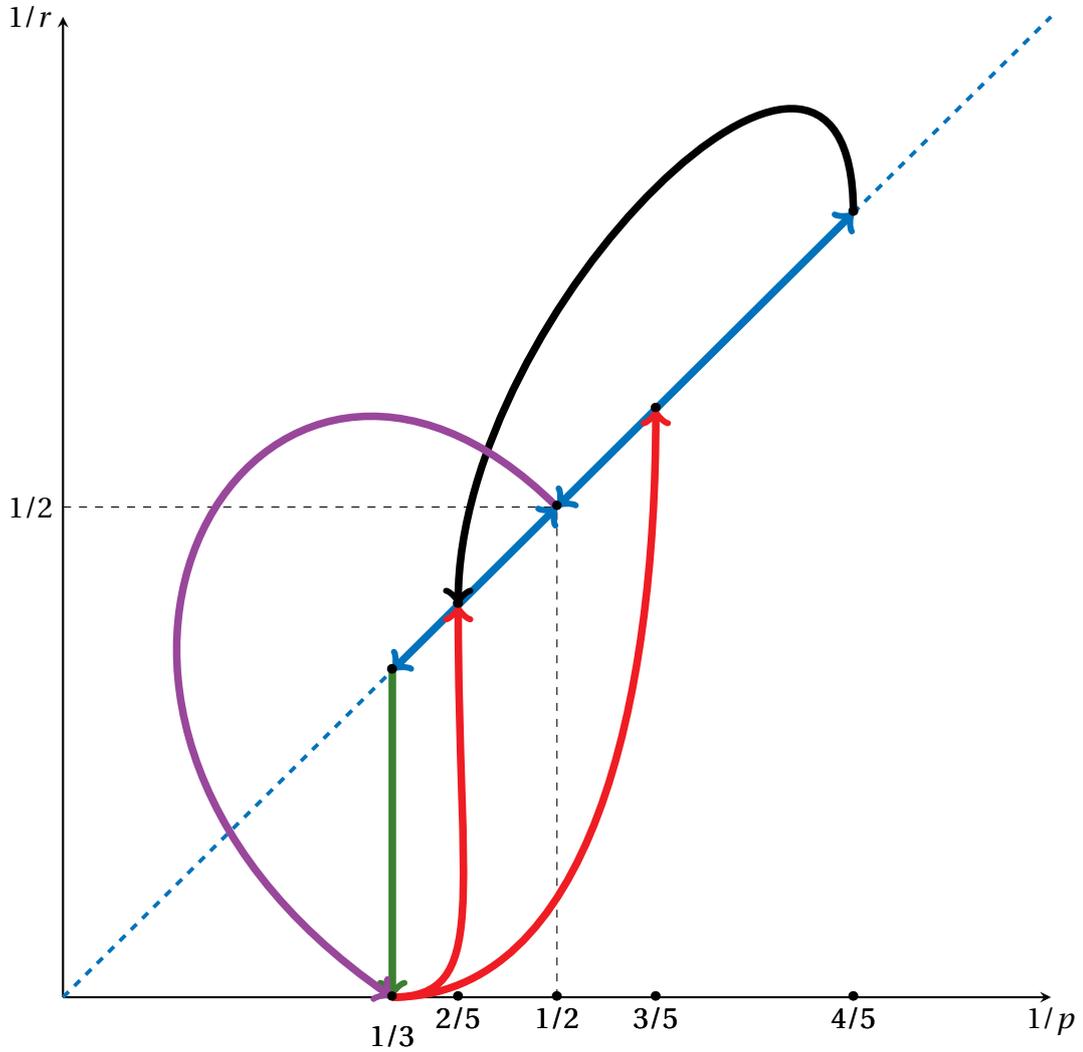

		The proof of  Theorem \ref{maintheorem} uses four important ingredients: a regularity result for  Burgers equation (Proposition \ref{thm1}), a higher regularity estimate for the  capillary Burgers equation (Proposition \ref{thm2}), 
		an energy estimate (Lemma \ref{Lem:ener}), and a result which allows us to ``increase'' the $r$ index of Besov spaces (Proposition \ref{thm3}). 
		Let us now sketch the proof, discarding lower-order terms (in particular all the terms containing $g=-|\partial_x|^{-1}\partial_x^2 u$) and taking borderline exponents in the estimates\footnote{Let us stress that we cannot reach these exponents since some of the constants (in particular the one in \eqref{explodeconst}) explode.}.
		The strategy is graphically represented in  Figure \ref{bellefigure}. 
		The starting point is Proposition \ref{thm1}, which for  $s=1$, $p=5/2$, $r=5/2$ and $\ksi=-|\partial_x|^{-1}\partial_x^4u$ (recall also \eqref{transforme g ksi}), roughly says that 
		\begin{equation*}
			\left\|u \right\|_{\B^{1/3}_{3,\infty}}^3\les \left\|u\right\|_{\B^1_{5/2,5/2}} \left\| u \right\|_{\B^3_{5/3,5/3}}. \label{111}
		\end{equation*}
		Using then the  interpolation inequality \eqref{Interpolation}, we get 
		\begin{equation*}
			\left\|u \right\|_{\B^{3}_{5/3,5/3}}\les \left\|u\right\|_{\B^{5}_{5/4,5/4}}^{1/3} \left\|u \right\|_{\B^2_{2,2}}^{2/3}.\label{222}
		\end{equation*}
		Proposition \ref{thm2} for $\alpha=2$, $p=5/4$, $q=5/2$ and therefore $\alpha'=1$, indicates that
		\begin{equation}\label{explodeconst}
			\left\|u\right\|_{\B^{5}_{5/4,5/4}} \les\left\|u \right\|_{\B^1_{5/2,5/2}}^{2}. 
		\end{equation}
		Using the interpolation inequality \eqref{Interpolation} once again, we find
		\begin{equation}\label{nonbuckle}
			\left\|u\right\|_{\B^1_{5/2,5/2}} \les \left\|u\right\|_{\B^{1/3}_{3,3}}^{3/5} \left\|u\right\|_{\B^2_{2,2}}^{2/5}. 
		\end{equation}
		From Lemma \ref{Lem:ener}, we obtain
		\begin{equation*}
			\left\|u\right\|_{\B^2_{2,2}}\les \left\|u\right\|_{\B^{1/3}_{3,\infty}}. \label{555}
		\end{equation*}

At this point, we see that we could have buckled the estimates if in  \eqref{nonbuckle}, the Besov norm $\left\|u\right\|_{\B^{1/3}_{3,3}}$ was replaced by the stronger norm $\left\|u\right\|_{\B^{1/3}_{3,\infty}}$. Unfortunately, this seems not doable with our method of proof.
		Therefore, we need Proposition \ref{thm3}   in order to control    $\left\|u\right\|_{\B^{1/3}_{3,3}}$ by $\left\|u\right\|_{\B^{1/3}_{3,\infty}}$. It is at this last stage that we lose a logarithm since \eqref{logarithm} gives
		\begin{equation*}
			\left\|u\right\|_{\B^{1/3}_{3,3}}\les \ln^{1/3}(L) \left\|u\right\|_{\B^{1/3}_{3,\infty}}.
		\end{equation*}
		 Putting all these estimates together, we find
		\begin{equation*}
			\left\|u \right\|_{\B^{1/3}_{3,\infty}} \les \ln^{1/2}(L),
		\end{equation*}
		which is \eqref{Main}.	\\
			
	As  mentioned, the first ingredient is an estimate for the inhomogeneous Burgers equation.  A similar estimate was obtained in \cite[Prop. 1]{Otto}. A related inequality  for the homogeneous Burgers equation has been recently derived  in \cite{Golse}.  

			Let us consider the following inhomogeneous Burgers equation:
			\begin{equation}\label{burger1}				
					\partial_t u + u \partial_x u = \left| \partial_x \right| g + \left| \partial_x \right| \ksi.
			\end{equation}
			\begin{proposition}[Besov estimate for the inhomogeneous Burgers equation]\label{thm1}
				Let $\ksi$, $g$ be smooth $L$-periodic functions. Then, for any smooth $L$-periodic solution $u$ of \eqref{burger1},
				there holds: For $s \in ]0,1[$, $r, r',p,p' \in [1, +\infty]$ verifying $\frac{1}{r}+\frac{1}{r'}=1$ and $\frac{1}{p}+\frac{1}{p'}=1$,  there exists a constant $c>0$ just depending on $s,r,p$ such that:
				\begin{align}\label{golse}
					\left\| u \right\|_{B^{1/3}_{3,\infty}}^3  \leq&
								c \left(\left\|u\right\|_{B^{1/3}_{3,\infty}} \left\| g \right\|_{B^{2/3}_{3/2,1}} 
								+ \left\| u \right\|_{B^s_{p,r}} \left\| \ksi \right\|_{B^{1-s}_{p',r'}} +\left\| u(0,\cdot) \right\|_{L^2[0,L]}^2  \right).
				\end{align}
				Therefore, taking the time-space average, it holds:
				\begin{align}\label{golse2}
					\left\| u \right\|_{\B^{1/3}_{3,\infty}}^3  \leq&
								c \left(\left\|u\right\|_{\B^{1/3}_{3,\infty}} \left\| g \right\|_{\B^{2/3}_{3/2,1}} 
								+ \left\| u \right\|_{\B^s_{p,r}} \left\| \ksi \right\|_{\B^{1-s}_{p',r'}}\right).
				\end{align}
				
			\end{proposition}
			The proof of \eqref{golse} is based on a modified  K\'arm\'an-Howarth-Monin identity:
		\begin{lemma}\label{monin}
		Let $\eta$ be a smooth $L$-periodic function and let $u$ be a smooth $L-$periodic solution with zero average of 
		\begin{equation}\label{burgeretastrat}\partial_t u + u\partial_x u = \eta,\end{equation}
		then for $h\in \R$,
		\begin{equation}\label{KHMstrat}
		\partial_t  \left(\frac{1}{2} \int_0^L |D^h u|D^hu dx\right)+\partial_h \frac{1}{6} \int_0^L |D^h u|^3 dx=\int_0^L D^h\eta  |D^hu| dx. 
		\end{equation}

		\end{lemma}
		  The usual K\'arm\'an-Howarth-Monin identity \cite{Frisch} states that 
		\begin{equation}\label{KHMoriginal}	\partial_t  \left(\frac{1}{2} \int_0^L \left(D^h u\right)^2 dx\right)+\partial_h \frac{1}{6} \int_0^L \left(D^h u\right)^3 dx=\int_0^L D^h\eta  D^h u dx. \end{equation}
		This formula can be easily checked by using  equation \eqref{burgeretastrat} and the periodicity. The main difference between \eqref{KHMoriginal} and \eqref{KHMstrat} is that in the latter, the coercive term $ \int_0^L |D^h u|^3 dx$ replaces the non-coercive term $\int_0^L \left(D^h u\right)^3 dx$. \\
We will give two proofs of \eqref{KHMstrat}. The first is by a direct computation and the second uses the kinetic formulation of Burgers equation following ideas of \cite{Golse}. Therefore, this second proof gives a new, and hopefully interesting, interpretation of the arguments of \cite{Golse}.

		The second ingredient is a higher regularity result for the capillary Burgers equation (see \cite[Prop. 2, p. 14]{Otto}).		
			\begin{proposition}[Higher regularity]\label{thm2}
				Let $p, q \in [1,+\infty[$, $\alpha \in \R$ satisfying:
				\begin{align*}
				 p+1\leq q \leq 2p,&&\text{ and }&& \alpha'=(6+\alpha)p/q-3 \in ]0,1[. 
				\end{align*}
				Then, there exists $c>0$ such that, if $u$, $g$ are smooth, $L$-periodic in $x$ and satisfy
				\begin{equation*}
					\partial_t u+ u\partial_x u+\partial_x^4u=\left|\partial_x \right|g,
				\end{equation*}
				the following estimate holds:
				\begin{equation}\label{higherr}
					\left\| u \right\|_{\B^{3+\alpha}_{p,p}}\leq c \left( \left\| u \right\|_{\B^{\alpha'}_{q,q}}^{q/p} +\left\| g \right\|_{\B^{\alpha}_{p,p}} \right).
				\end{equation}
			\end{proposition}
			Proposition \ref{thm2}  allows to jump from higher derivatives to smaller ones in Besov spaces. The proof, which we will not provide, is based on a narrow-band Littlewood-Paley decomposition.
			
			The third ingredient is an elementary energy estimate, which directly bounds the $L^2\left([0,T],H^2([0,L])\right) \cong B^2_{2,2}$ norm of a solution $u$ of the inhomogeneous capillary Burgers equation.
			\begin{lemma}[Energy estimate]\label{Lem:ener}
				Let $u$ be a smooth solution of:
				\begin{equation}\label{der4}
					\partial_t u + u\partial_x u +\partial_x^4 u = \left| \partial_x \right| g.
				\end{equation}
				 Then, the following estimate holds:
				\begin{equation}\label{JJ}
					\left\|u\right\|_{\B^2_{2,2}} \leq c \left\|u \right\|_{\B^{1/3}_{3,\infty}}^{1/2} \left\| g\right\|_{\B^{2/3}_{3/2,1}}^{1/2}.
				\end{equation}
			\end{lemma}
			\begin{proof} Since the proof is straightforward, we give it now. Integrating over  the equation \eqref{der4} over $[0,T]\times[0,L]$, we get:
				\begin{equation*}
					\int_0^L u(T,x)^2dx - \int_0^L u(0,x)^2dx =\int_0^T\int_0^L u\left|\partial_x\right|gdxdt-\int_0^T\int_0^L(\partial_x^2 u)^2dxdt.
				\end{equation*}
				Therefore by \eqref{Besov1}
				\begin{equation*}
					\int_0^T\int_0^L\left(\partial_x^2u \right)^2dxdt \leq \int_0^Lu(0,x)^2 dx+ \left\|u\right\|_{B^{1/3}_{3,\infty}}\left\|g\right\|_{B^{2/3}_{3/2,1}},
				\end{equation*}
				 taking then the time-space average yields the result.
		\end{proof}
		
			As already mentioned, these three estimates will not be sufficient to conclude. We will also need an estimate relating Besov norms with different exponents $r$. One can easily see that  if $r_1>r_2$, then $\left\|\cdot \right\|_{\B^s_{p,r_1}} \leq \left\|\cdot \right\|_{\B^{s}_{p,r_2}}$ (it is a consequence of convexity inequality). In fact, it is possible to reverse the inequality for {\it solutions} of \eqref{Kuramoto}, but this comes with a price: a logarithm of the spatial period $L$ appears.
			\begin{proposition}[Increasing the index $r$] \label{thm3}
				There exists $c >0$ such that, for all $L\ge 2$, $u$ solution of \eqref{Kuramoto}, the following estimate holds:
				\begin{align}\label{logarithm}
					\left\|u \right\|_{\B^{1/3}_{3,3}} \leq c \ln^{1/3}(L)\left\|u\right\|_{\B^{1/3}_{3,\infty}}.
				\end{align}
			\end{proposition}

\section{Proof of Theorem \ref{maintheorem}}\label{sec:proofmain}
	In this section, we derive the main theorem from the above propositions. We now consider the rescaled Besov norm $\B^s_{p,r}$ as a point $(s,1/p,1/r)$ in the space $\R^3$. All the norms involved in our problem lie in the rectangle $\mathcal{P}$ of $\R^3$ defined by\footnote{See Figure \ref{bellefigure} in Section 3.2 which represents the strategy in $\mathcal{P}$.}:
				\begin{equation*}
					 \left\{ \begin{array}{r c l}
						s&=& 10/p-3, \\
						\frac{1}{p}&\in& [0,1],\\
						 \frac{1}{r}& \in & [0,1].
					\end{array}
					\right.
				\end{equation*}
				\begin{proof}[Proof of Theorem \ref{maintheorem}]
	Let $u$ be a solution of \eqref{Kuramoto}.
	It is convenient to introduce the abbreviation:
	\begin{equation*}
		D(\alpha)= \left\| u \right\|_{\B^{10\alpha-3}_{\alpha^{-1},\alpha^{-1}}} \qquad \Dast=\left\| u \right\|_{\B^{1/3}_{3,\infty}}.
	\end{equation*}
	Notice that $\Dtw=\left\| u \right\|_{\B^{1/2}_{2,2}}$ and $\Dth=\left\| u \right\|_{\B^{1/3}_{3,3}}$. With this notation, interpolation inequality \eqref{Interpolation} takes the form
	\begin{equation}\label{interpolD}
	D(\alpha)\le D^{\theta}(\alpha_1)D^{1-\theta}(\alpha_2)
	\end{equation}
	for $\alpha=\theta \alpha_1 + (1-\theta)\alpha_2$ and $\theta\in[0,1]$.
	
	Letting $s=10 \alpha-3$, $p=\alpha^{-1}$ and $r=\alpha^{-1}$ in \eqref{golse2} and using \eqref{transforme g ksi} for $\xi=-|\partial_x|^{-1} \partial_x^4 u$, \eqref{golse2} can   be rewritten as
	\begin{equation}\label{golserew}
	\Dastp{3}\le c \left(D(\alpha)D(1-\alpha)+ \Dast \| g \|_{\B^{2/3}_{3/2,1}}\right),
	\end{equation}
	for  $\alpha \in ]3/10,2/5[$. In turn, \eqref{higherr} with $p^{-1}=\beta$, $q^{-1}=\gamma$ and $\alpha=10\beta-6$ (which implies $\alpha'=10\gamma-3$)  gives
	\begin{equation}\label{highrew}
	D(\beta)\le c \left(D^{\beta/\gamma}(\gamma)+ \|g\|_{\B^{10\beta-6}_{\beta^{-1},\beta^{-1}}}\right),
	\end{equation}
	 for    $\gamma \in ]3/10, 2/5[$ and $\frac{\gamma}{1-\gamma}\le \beta\le 2\gamma$. Finally, \eqref{JJ}  is equivalent to
	\begin{equation}\label{JJrew}
	\Dtw\le c \Dastp{1/2}\|g \|_{\B^{2/3}_{3/2,1}}^{1/2},
	\end{equation}
	and \eqref{logarithm} to
	\begin{equation}\label{logarithmrew}
	\Dth\le c \ln^{1/3}(L) \Dast.
	\end{equation}
	
	Our first goal is to argue that for $g=- |\partial_x|^{-1} \partial_x^2 u$, we can replace in the above estimates all the Besov norms involving $g$ by $\Dast$. By \eqref{transforme g ksi} and \eqref{Interpolation},
	\begin{equation*}\|g\|_{\B^{10\beta-6}_{\beta^{-1},\beta^{-1}}}\le c\|u\|_{\B^{10\beta-5}_{\beta^{-1},\beta^{-1}}}
	\le c \|u\|_{\B^{10\beta-3}_{\beta^{-1},\beta^{-1}}}^{1/2}\|u\|_{\B^{10\beta-7}_{\beta^{-1},\beta^{-1}}}^{1/2}=c D^{1/2}(\beta)\|u\|_{\B^{10\beta-7}_{\beta^{-1},\beta^{-1}}}^{1/2}.
	\end{equation*}
	Hence, in view of \eqref{JJrew}, Young's inequality and since by \eqref{transforme g ksi}, $\|g \|_{\B^{2/3}_{3/2,1}}\le c\|u \|_{\B^{5/3}_{3/2,1}}$,  it will be enough to prove that 
	\begin{equation}\label{removeg}  \|u\|_{\B^{10\beta-7}_{\beta^{-1},\beta^{-1}}}+\|u \|_{\B^{5/3}_{3/2,1}}\le c (\Dast+\Dtw)\end{equation}
for $\beta\in]11/15,4/5[$ (which reduces the use of \eqref{highrew} to $\gamma\in]1/3, 2/5[$). We can indeed prove more generally that for $1/3<s<2$, $p\le 2$ and any $q\ge 1$, there holds
\begin{equation}\label{interpolproof}
\|u\|_{\B^s_{p,q}}\le c \left( \|u\|_{\B^{2}_{2,2}}+ \|u\|_{\B^{1/3}_{3,\infty}}\right).
\end{equation}	
	Thanks to Jensen's inequality, we have $\|u\|_{\B^s_{p,q}}\le \|u\|_{\B^s_{2,q}}$ and $\|u\|_{\B^{1/3}_{3,\infty}}\ge \|u\|_{\B^{1/3}_{2,\infty}}$. 
	By monotonicity of the Besov norms with respect to the last index, there also holds $\|u\|_{\B^s_{2,q}}\le \|u\|_{\B^s_{2,1}}$ and $ \|u\|_{\B^2_{2,2}}\ge \|u\|_{\B^2_{2,\infty}} $. Therefore, we are left with proving that 
	\[\|u\|_{\B^s_{2,1}}\le c\left(\|u\|_{\B^{1/3}_{2,\infty}}+ \|u\|_{\B^2_{2,\infty}}\right).\]
	By definition of the Besov norms, for $1/3<s<2$,
	\begin{align*}
	\|u\|_{B^s_{2,1}}&=\sum_{k\ge 0} 2^{k}\|u_k\|_{L^2}+\sum_{k< 0} 2^{k}\|u_k\|_{L^2} \\
	&\le \sup_{k} 2^{2k} \|u_k\|_{L^2}\sum_{k\ge 0} 2^{-k} + \sup_{k} 2^{\frac{1}{3}k} \|u_k\|_{L^2}\sum_{k\le 0} 2^{\frac{2}{3}k}\\
	&\le c \left(\|u\|_{B^{1/3}_{2,\infty}}+ \|u\|_{B^2_{2,\infty}}\right),
\end{align*}
which after taking the average over time and space, finishes the proof of \eqref{interpolproof}. \\

To sum up, we now have that \eqref{golserew}, \eqref{highrew} and \eqref{JJrew} together with \eqref{removeg} imply

		\begin{equation}\label{golserew2}
	\Dastp{3}\le c \left(D(\alpha)D(1-\alpha)+ \Dastp{2}\right)
	\end{equation}
	for  $\alpha \in ]3/10,2/5[$, 
	\begin{equation}\label{highrew2}
	D(\beta)\le c \left(D^{\beta/\gamma}(\gamma)+ \Dast\right)
	\end{equation}
	 for    $\gamma \in ]1/3, 2/5[$ and $\frac{\gamma}{1-\gamma}\le \beta\le 2\gamma$ and
	\begin{equation}\label{JJrew2}
	\Dtw\le c \Dast. 
	\end{equation}

We now gather the above estimates in order to bound $\Dast$. Passing to the logarithm in the above inequalities, we see that optimizing the parameters to get the best power of $\ln(L)$ is equivalent to a linear programming problem. Its solution thus lie at the boundaries of the admissible domain. It is not hard to see that in particular, we want to take $\frac{\beta}{2}=\gamma=\alpha$   with $\alpha$ as close as possible to $2/5$. Let $\theta,\eta\in]0,1[$ be such that 
\[ (1-\alpha)= \theta \beta + (1-\theta) \frac{1}{2} \qquad \textrm{and } \qquad  \alpha=\eta \frac{1}{2}+ (1-\eta) \frac{1}{3},\]
so that $\theta$ is close to $1/3$ and $\eta$ is close to $2/5$. Thanks to \eqref{interpolD},
\begin{equation}\label{estiminterpolmain}
D(1-\alpha)\le  D^{\theta}(\beta) D^{1-\theta}(1/2) \qquad \textrm{and } \qquad D(\alpha)\le  D^{\eta}(1/2) D^{1-\eta}(1/3).
\end{equation}  
		Since we can assume that $\Dast\ge 1$, we get from \eqref{golserew2}, \eqref{estiminterpolmain} and \eqref{highrew2},
		\begin{align*}
		\Dastp{3}&\le c D(\alpha) \left(D^{2\theta}(\alpha)+\Dastp{\theta}\right) D^{1-\theta}(1/2)\\
		&\le c \left( D^{1+2\theta}(\alpha) \Dastp{1-\theta} + D(\alpha)\Dast \right),
		\end{align*}
		where in the last inequality, we used \eqref{JJrew2}. From \eqref{estiminterpolmain}, \eqref{JJrew2} and \eqref{logarithmrew}, we deduce
		
		\[
		\Dastp{3}\le c\left( \Dastp{2+\theta}\ln^{\frac{1}{3}(1+2\theta)(1-\eta)}(L) + \Dastp{2} \ln^{\frac{1-\eta}{3}}(L)\right).
		\]
		Dividing by $\Dastp{2}$ this inequality and noticing that for $\eta$ close to $2/5$, $\frac{1-\eta}{3}$ is close to $1/5$, we obtain that if $\Dast\ge \ln^{\frac{1-\eta}{3}}(L)$, then 
		\[\Dast\le c\Dastp{\theta} \ln^{\frac{1}{3}(1+2\theta)(1-\eta)}(L), \]
		which gives finally
		\[\Dast\le c \ln^{\frac{1}{3}\frac{(1+2\theta)(1-\eta)}{1-\theta}}(L) \]
and thus the result since $\lim_{\theta\uparrow 1/3,\eta\uparrow 2/5} \ \frac{1}{3}\frac{(1+2\theta)(1-\eta)}{1-\theta}=1/2$.
	\end{proof}
\section{Proof of Proposition \ref{thm1}}\label{sec:proofburgers}
	For the reader's convenience, let us recall the statement of Proposition \ref{thm1}. Let $u$ be a smooth solution of the
 following inhomogeneous Burgers equation:
			\begin{equation}\label{burger2}
					\partial_t u + u \partial_x u = \left| \partial_x \right| g + \left| \partial_x \right| \ksi.
			\end{equation}
			
			\begin{proposition*}
				Let $\ksi$, $g$ be smooth $L$-periodic functions. Then, for any smooth $L$-periodic solution $u$ of \eqref{burger2},
				there holds: For $s \in ]0,1[$, $r, r', p, p' \in [1, +\infty]$ verifying $\frac{1}{r}+\frac{1}{r'}=1$ and $\frac{1}{p}+\frac{1}{p'}=1$,  there exists a constant $c \in \R_+^*$ just depending on $s,r,p$ such that:
				\begin{align*}
					\left\| u \right\|_{B^{1/3}_{3,\infty}}^3  \leq&
								c \left(\left\|u\right\|_{B^{1/3}_{3,\infty}} \left\| g \right\|_{B^{2/3}_{3/2,1}} 
								+ \left\| u \right\|_{B^s_{p,r}} \left\| \ksi \right\|_{B^{1-s}_{p',r'}} +\left\| u(0,\cdot) \right\|_{L^2[0,L]}^2  \right).
				\end{align*}
			\end{proposition*}
		Before proceeding further, let us remark that, by approximation, this applies to any (possibly non smooth) entropy solution of Burgers equation  \eqref{burger2}. Indeed, if we consider a solution $u$ of 
		\[
		\partial_t u - u\partial_x u -\eps \partial^2_x u=0, 
		\]
		then it is a smooth solution of \eqref{burger2}  with $g=0$ and $\xi= \eps \left|\partial_x\right| u$.  A careful inspection of the proof of Proposition \ref{thm1} shows that for $p=r=2$ it extends to $s=1$, yielding
		\[
		 \left\| u \right\|_{B^{1/3}_{3,\infty}}^3  \leq c \left( \left\| u \right\|_{B^{1}_{2,2}}\left\| \ksi \right\|_{B^{0}_{2,2}} +\left\| u(0,\cdot) \right\|_{L^2[0,L]}^2 \right),
		\]
that is 
\[ \left\| u \right\|_{B^{1/3}_{3,\infty}}^3  \leq c \left( \eps \left\| \partial_x u \right\|^2_{L^2} +\left\| u(0,\cdot) \right\|_{L^2[0,L]}^2 \right).
		\]
		Combining this with the energy inequality: $\eps  \left\| \partial_x u \right\|^2_{L^2} \le c \left\| u(0,\cdot) \right\|_{L^2[0,L]}^2$, gives
		\[ \left\| u \right\|_{B^{1/3}_{3,\infty}}^3  \leq c \left\| u(0,\cdot) \right\|_{L^2[0,L]}^2,\]
		which passes to the limit as $\eps\to 0$.
The indices are optimal in the light of the result of De Lellis and Westdickenberg  \cite{DeLellis} which states that we cannot hope to have more regularity, in the sense that the Besov index $s$ cannot be better than $1/3$.\\
		
		As already pointed out, the proof of the aimed estimate is based on a modified K\'arm\'an-Howarth-Monin identity:
		\begin{lemma}
		Let $\eta$ be a smooth $L$-periodic function and let $u$ be a smooth $L-$periodic solution with zero average of 
		\begin{equation}\label{burgereta}\partial_t u + u\partial_x u = \eta\end{equation}
		then for $h\in \R$,
		\begin{equation}\label{KHM}
		\partial_t  \left(\frac{1}{2} \int_0^L|D^h u|D^hu dx\right)+\partial_h \frac{1}{6} \int_0^L \left|D^h u\right|^3 dx=\int_0^L D^h\eta  \left|D^hu\right| dx. 
		\end{equation}
	
		\end{lemma}
		\begin{proof}
		By periodicity,  \eqref{KHM} will be a direct	 consequence of the following pointwise identity:
			\begin{equation}\label{pointwiseeq}
				\frac{1}{2}\partial_t\left(|D^h u|D^hu \right)+ \frac{1}{6} \partial_h \left|D^h u\right|^3 +\frac{1}{2} \partial_x\left( u |D^h u| D^h u+\frac{1}{3} |D^h u|^3\right) =D^h\eta  \left|D^hu\right|.
			\end{equation}
			For simplicity, let us introduce the notation $u^h(x)=u(x+h)$ (so that $D^h u=u^h-u$). Using \eqref{burgereta} we get:
			\begin{align*}
				\frac{1}{2}\partial_t\left( |D^h u|D^h u\right)+ \frac{1}{6} \partial_h \left|D^h u\right|^3&=|D^h u|\partial_t (D^h u)+\frac{1}{2}\left|D^hu \right|D^hu\partial_xu^h\\
				&=\left|D^hu\right|D^h\left(\eta-u\partial_xu\right)+\frac{1}{2}\left|D^hu\right|\left(u^h\partial_xu^h-u\partial_xu^h\right)\\
				&=D^h\eta  \left|D^hu\right|+\frac{1}{2}\left|D^hu\right|\left(2u\partial_xu-u\partial_xu^h-u^ h\partial_xu^h\right).
			\end{align*}
			It remains to prove that
			\begin{equation}\label{equaldx}
			|D^hu|\left(2u\partial_xu-u\partial_xu^h-u^ h\partial_xu^h\right)=-\partial_x\left( u |D^h u| D^h u+\frac{1}{3} |D^h u|^3\right).\end{equation}
			We start with 
			\[|D^hu|2u\partial_xu= -u^2 \partial_x |D^h u| + \partial_x \left( |D^h u| u^2 \right)\]
			and
			\begin{align*}|D^hu|u \partial_xu^h=&-\partial_x \left(|D^h u| u\right) u^h+ \partial_x \left(|D^h u| u u^h\right)\\
			=&-u u^h\partial_x |D^h u| -u^h|D^hu|\partial_x u+ \partial_x \left(|D^h u| u u^h\right), 
			\end{align*}
			to get
			\begin{align*}
			|D^hu|\left(2u\partial_xu-u\partial_xu^h-u^ h\partial_xu^h\right)=&
			-u^2 \partial_x |D^h u| + \partial_x \left( |D^h u| u^2 \right)+u u^h\partial_x |D^h u| \\
			& +u^h|D^hu|\partial_x u- \partial_x \left(|D^h u| u u^h\right) -|D^h u|u^h \partial_x u^h\\
			=&-\partial_x \left( u|D^h u| D^h u\right)+u D^h u\partial_x |D^h u|- u^h |D^h u| \partial_x D^h u. 
			\end{align*}
			But since
			\[
				D^hu\partial_x |D^hu |=|D^hu |\partial_x D^hu,
			\]
			then
			 \begin{align*}
			|D^hu|\left(2u\partial_xu-u\partial_xu^h-u^ h\partial_xu^h\right)&=-\partial_x \left( u |D^h u|D^h u \right)-  |D^h u|D^h u \partial_x D^h u\\
			&=-\partial_x\left( u |D^h u| D^h u+\frac{1}{3} |D^h u|^3\right),
			\end{align*}
			which concludes the proof of \eqref{equaldx}.
		\end{proof}
		
		\begin{remark}
		 Arguing along the same lines, one can prove that more generally, if $a$ is non-negative and if $u$ is a smooth solution with zero average of 
		 \begin{equation}\label{conserva}\partial_t u+\partial_x [a(u)]= \eta\end{equation}
		 then
		 \begin{multline*}\partial_t  \left(\frac{1}{2} \int_0^L|D^h u|D^hu dx\right)+\partial_h  \left(\int_0^L |D^hu| (a(u)+a(u^h))-2|A(u^h)-A(u)| dx\right)\\
		 =\int_0^L D^h\eta  \left|D^hu\right| dx, 
		 \end{multline*}
where $A'=a$. Notice that Burgers equation corresponds to \eqref{conserva} with $a(u)=\frac{1}{2} u^2$. If $a$ is $C^1$ and monotone in the sense that there exist $\beta\ge 1$ and $C>0$, such that for $ v\ge w$,
\[a'(v)-a'(w)\ge  C (v-w)^\beta,\]
then one can obtain a similar estimate to \eqref{golse} by using that for $\bar{u}\ge u$
\begin{align*} (\bar{u}-u)(a(\bar{u})+a(u))-2(A(\bar{u})-A(u))&=\int_u^{\bar{u}} \int_{w}^{\bar{u}} (a'(v)-a'(w)) dv dw\\
&\ge C |\bar{u}-u|^{\beta+2}.
\end{align*}
In this way, one can fully recover the results from \cite{Golse}.
		\end{remark}

		We now give an alternative proof of (the integrated form of) \eqref{KHM}  following ideas from \cite{Golse}. 
		This proof  uses the kinetic formulation of the inhomogeneous Burgers equation together with the use of an interaction identity (see \eqref{InteractionId}).

\begin{lemma}
	Let $\eta$ be a smooth $L$-periodic function and let $u$ be a smooth $L-$periodic solution with zero average of 
		\begin{equation}\label{burgereta2}\partial_t u + u\partial_x u = \eta,\end{equation}
		then for $h\in \R$,
		\begin{multline}\label{KHMint}
	\left[\frac{1}{2} \int_0^h \int_0^L |D^\Delta u|D^\Delta u dx d\Delta \right]_0^T+ \frac{1}{6} \int_0^T\int_0^L \left|D^h u\right|^3 dx dt \\=\int_0^T\int_0^L\int_0^h D^\Delta \eta   \left|D^\Delta u\right|  d\Delta dx dt.\end{multline} 
		\end{lemma}
\begin{proof}
Before starting the proof, let us point out that since we have a direct proof of \eqref{KHM}, we will not take care of regularity issues. Nevertheless, all passages can be rigorously justified by a suitable approximation argument.\\
\medskip

{\it Step 1.} Without loss of generality, we can assume that $h>0$. Letting:
			\begin{equation*}
				f(t,x,v)=\left\{
					\begin{array}{r l}
						1	& \text{if }	v \leq u(t,x), \\
						0	& \text{if } v> u(t,x),\\				
					\end{array}
				\right.
			\end{equation*}
			equation \eqref{burgereta2} is equivalent to the following kinetic formulation (see \cite{Lions}):
			\begin{equation}\label{kinetic}
				\partial_t f(t,x,v) + v \partial_x f(t,x,v) = - \partial_v f(t,x,v) \eta (t,x).
			\end{equation}
Notice that  since $u$ is bounded,  $D^h f$ is integrable even though $f$ is not. We are going to compute only integrals depending on $D^h f$ and will therefore not have to deal with integrability issues.
		As in \cite[Lem. 4.3]{Golse}, we define:
			\begin{equation}\nonumber
				M_u(v)= \left\{
					\begin{array}{r l}
						1	& \text{if }	v \leq u, \\
						0	& \text{if } v> u.\\				
					\end{array}
				\right.
			\end{equation}
			
			We first claim that for all $u$, $\bar{u} \in \R$ the following equality holds:
				\begin{equation}\label{estimate}
					\frac{1}{6}|u-\bar{u}|^3 = \int_{\R} \int_{\R} \left[\mathbb{1}_{\R^+}(v-w) \right](v-w) \left(M_u(v)-M_{\bar{u}}(v)\right)\left(M_u(w)-M_{\bar{u}}(w)\right)dv dw.
				\end{equation}
				Without loss of generality, we can suppose that $u \geq \bar{u}$. Then:
				\begin{equation*}
					M_u(v)-M_{\bar{u}}(v)=\mathbb{1}_{]\bar{u},u[}(v). 
				\end{equation*}
				Thus \eqref{estimate} follows from:
\begin{align*}
					\int_{\R}\int_{\R} \left[\mathbb{1}_{\R_+}(v-w) \right](v&-w) \left(M_u(v)-M_{\bar{u}}(v) \right)\left( M_{u}(w)-M_{\bar{u}}(w) \right)dvdw\\
					&=\int_{\R}\int_{\R} \left[\mathbb{1}_{\R_+}(v-w)\mathbb{1}_{]\bar u,u[}(v)\mathbb{1}_{]\bar u,u[}(w)\right] (v-w)dvdw\\
					&=\int_{\bar u}^u\int_w^u (v-w) dvdw=\frac{1}{6} \left|u-\bar u\right|^3.
				\end{align*}

			Letting
			\begin{equation*}
				Q(h)=\int_0^T \int_0^L \int_{\R\times\R} \left[\mathbb{1}_{\R_+}(v-w)\right](v-w)D^h f(t,x,v)D^h f(t,x,w)dvdwdxdt,
			\end{equation*}
			we see that proving \eqref{KHMint} is equivalent to 
			\begin{equation}\label{KHMintequiv}
			Q(h)=\int_0^T\int_0^L\int_0^h D^\Delta \eta  \left|D^\Delta u\right| d\Delta dx dt -\left[\frac{1}{2} \int_0^h \int_0^L |D^\Delta u|D^\Delta u dx d\Delta \right]_0^T. 
			\end{equation}
			\medskip

	{\it Step 2.}		To cope with the quantity $Q$, the main tool is the following interaction identity (see \cite{Golse}), which have been introduced first by Varadhan (\cite[Lem. 22.1]{Tartar}):
			Let $A$, $B$, $C$, $D$, $E$, $F$ $:[0,T]\times [0,L] \rightarrow \R$ be functions satisfying the following system
			\begin{equation}\label{interassumption1}
			\left\{
				\begin{array}{c}
					\partial_t A +\partial_x B = C, \\
					\partial_t D + \partial_x E = F,
				\end{array}
			\right.
			\end{equation}
			and having zero spatial average. Then the following identity holds:
			\begin{align}\label{InteractionId}
				\int_0^T \int_0^L (AE-BD)=&\int_0^T\int_0^LA(t,x) \left(\int_0^x F(t,y)dy \right)dxdt \nonumber \\
					&+ \int_0^T\int_0^LC(t,x)\left(\int_0^xD(t,y)dy\right)dxdt \\
					&-\left[\int_0^L\int_0^x A(t,x)D(t,y)dydx \right]_{t=0}^{t=T}.\nonumber
			\end{align}
			
				Indeed, by Taylor expansion:
				\begin{align*}
					\int_0^T \int_0^L A(t,x)E(t,x) dx dt =
						&\int_0^T \int_0^L \int_0^x A(t,x)\partial_xE(t,y)dydxdt\\
						&-\int_0^T\int_0^L A(t,x)E(t,0)dxdt.
				\end{align*}
				Since $A$ has zero spatial average, the second term vanishes. Using equation \eqref{interassumption1} to compute the first term and integrating by parts, we get:
				\begin{align*}
					\int_0^T\int_0^L\int_0^xA(t,x)\partial_xE(t,y)dydxdt =& \int_0^T\int_0^L\int_0^xA(t,x)F(t,y)dydxdt\\
					&-\int_0^T\int_0^L\int_0^x A(t,x)\partial_tD(t,y)dydxdt \\
					=&\int_0^T\int_0^L\int_0^x A(t,x)F(t,y)dydxdt\\
					&+\int_0^T\int_0^L\int_0^x \partial_tA(t,x)D(t,y)dydxdt\\
					&-\left[\int_0^L\int_0^x A(t,x)D(t,y)dydx \right]_{t=0}^{t=T}.
				\end{align*}
				Let us now compute more precisely the second term, using \eqref{interassumption1}:
				\begin{align*}
					\int_0^T\int_0^L\int_0^x \partial_tA(t,x)D(t,y)dydxdt=&\int_0^T\int_0^L\int_0^xC(t,x)D(t,y)dydxdt\\
					&-\int_0^T \int_0^L\partial_xB(t,x) \int_0^x D(t,y)dydxdt\\
					=&\int_0^T\int_0^L\int_0^xC(t,x)D(t,y)dydxdt\\
					&+\int_0^T\int_0^L B(t,x)D(t,x)dxdt,
				\end{align*}
				which concludes the proof of \eqref{InteractionId}.
			\medskip

			{\it Step 3.}
			We apply the interaction identity to:
			\begin{equation*}
				\left\{
				\begin{array}{l c l}
				A(t,x,v)=D^h  f(t,x,v), & & D(t,x,w)=A(t,x,w),\\
				B(t,x,v)=v D^h f(t,x,v), & & E(t,x,w)= B(t,x,w), \\
				C(t,x,v)=- \partial_v D^h   (f(t,x,v)\eta(t,x)), && F(t,x,w)=C(t,x,w).\\
				\end{array}
				\right.
			\end{equation*}
			 Note that $A$,$B$,$C$,$D$,$E$,$F$ implicitly depend on $h$. Multiplying each side of the identity by $\mathbb{1}_{\R_+}(v-w)$ and integrating it, we get:
			\begin{align*}
				Q(h)=&-\int_{\R\times\R} \mathbb{1}_{\R_+}(v-w) \int_0^T \int_0^L \left(A(t,x,v)E(t,x,w)-B(t,x,v)D(t,x,w) \right)dxdtdvdw\\
					=&-\int_{\R\times\R} \mathbb{1}_{\R_+}(v-w) \int_0^T \int_0^L A(t,x,v) \left(\int_0^x F(t,y,w)dy \right)dxdtdvdw \\
					&-\int_{\R\times\R} \mathbb{1}_{\R_+}(v-w) \int_0^T \int_0^L C(t,x,v) \left( \int_0^x D(t,y,w) dy \right)dxdtdvdw\\
					&+\int_{\R\times\R} \mathbb{1}_{\R_+}(v-w) \left[ \int_0^L\int_0^x A(t,x,v)D(t,y,w)dydx \right]_{t=0}^{t=T}dvdw.\\
					=&-\int_{\R\times\R} \mathbb{1}_{\R_+}(v-w) \int_0^T \int_0^L A(t,x,v) \left(\int_0^x F(t,y,w)dy \right)dxdtdvdw \\
					&-\int_{\R\times\R} \mathbb{1}_{\R_+}(w-v) \int_0^T \int_0^L F(t,y,w) \left( \int_0^y A(t,x,v) dx \right)dydtdvdw\\
					&+\int_{\R\times\R} \mathbb{1}_{\R_+}(v-w) \left[ \int_0^L\int_0^x A(t,x,v)D(t,y,w)dydx \right]_{t=0}^{t=T}dvdw.\\
					=&Q_1+Q_2+Q_3.
 					\end{align*}
But, by periodicity, $\int_0^y A(t,x,w) dx=-\int_y^{L} A(t,x,w) dx$, and thus
\begin{align*}
 \int_0^L F(t,y,v) \left( \int_0^y A(t,x,w) dx \right)dy&=-\int_0^L A(t,x,w)\left(\int_0^x F(t,y)d y\right) dx.
\end{align*}
Therefore,
\[
 Q_2=\int_{\R\times\R} \mathbb{1}_{\R_+}(w-v) \int_0^T \int_0^L A(t,x,v) \left(\int_0^x F(t,y,w)dy \right)dxdtdvdw.
\]
\medskip

{\it Step 4.}
In the next two steps,  the time variable plays no role. We will therefore consider 
			\[\overline{Q}_1=-\int_{\R\times\R} \mathbb{1}_{\R_+}(v-w) \int_0^L A(x,v) \left(\int_0^x C(y,w)dy \right)dxdvdw\]
and 
\[\overline{Q}_2=\int_{\R\times\R} \mathbb{1}_{\R_+}(w-v) \int_0^L A(x,v) \left(\int_0^x C(y,w)dy \right)dxdvdw.\]
By definition of $A$ and $C$, we have
\begin{align*}
 \overline{Q}_1-\overline{Q}_2&=-\int_{\R\times\R} \int_0^L D^hf(x,v) \left(\int_0^x \partial_v D^h(f(y,w) \eta(y))dy \right)dxdvdw\\
 &=\int_0^L D^h u(x) \left(\int_0^x D^h \eta(y) dy\right)dx.
\end{align*}
 The y-integral then telescopes to:
				\[
					 \overline{Q}_1-\overline{Q}_2=\int_0^L D^h u(x) \left(\int_x^{x+h}-\int_0^h\right)  \eta(y)dydx.
				\]
				Since by periodicity we have,
				\begin{equation*}
					\int_0^L D^h  u(x)dx=0,
				\end{equation*}
				this reduces to
				\begin{align}
					\overline{Q}_1-\overline{Q}_2&=\int_0^L D^h u(x) \int_x^{x+h} \eta(y)dydx\nonumber\\
					&=\int_0^L \eta(y)\int_{y-h}^y D^h u(x)dx dy \nonumber\\
					&=\int_0^L \int_0^h \eta(y) (u(y+\Delta)-u(y-\Delta)) dy d\Delta\label{Q1Q2}.
				\end{align}

				\medskip
				
		{\it Step 5.} Here we argue that 
		\begin{equation}\label{Q1}
		\overline{Q}_1+\overline{Q}_2=\int_0^L\int_0^h D^\Delta \eta  \left|D^\Delta u\right|  d\Delta dx.
		\end{equation}
			For this we prove first that 
			\begin{equation}\label{estimQ1}
			 \overline{Q}_1=\frac{1}{2}\int_0^h\int_0^L  D^{\Delta} \eta \left|D^{\Delta} u\right| +  \eta (u(x-\Delta)-u(x+\Delta)) dx d\Delta.
			\end{equation}
Indeed, combined with \eqref{Q1Q2}, this would give,
\[\overline{Q}_1+\overline{Q}_2= 2\overline{Q_1}-(\overline{Q}_1-\overline{Q}_2)= \int_0^L\int_0^h D^\Delta \eta  \left|D^\Delta u\right|  d\Delta dx.\]
which is \eqref{Q1}.
				By definition of $A$ and $C$:
				\begin{align*}
					\overline{Q}_1=&\int_{\R\times\R} \mathbb{1}_{\R_+}(v-w) \int_0^L D^h f(x,v) \left(\int_0^x \partial_v D^h \left(f(y,w)\eta(y)\right) dy\right)dxdvdw\\
					=&\int_{-\infty}^{+\infty}\int_0^L D^h f(x,v) \int_{-\infty}^v \int_0^x \partial_v D^h\left(f(y,w)\eta(y)\right)dydwdxdv\\
					=&\int_{-\infty}^{+\infty}\int_0^L D^h f(x,v) \int_0^x D^h\left(f(y,v)\eta(y)\right) -D^h \eta(y) dydxdv.
				\end{align*}
Arguing as in {\it Step 4}, we find
				\[
					\overline{Q}_1=  \int_0^L \eta(y) \left(\int_{0}^h \int_{-\infty}^{+\infty}D^h f(y-x,v)\left(f(y,v)-1\right)dxdv\right)dy.
				\]
				After changing the names of the variables to $y=\tilde{x}$ and $x=\Delta$ we obtain 
				\[\overline{Q}_1=  \int_0^L \eta(\tilde x) \left(\int_{0}^h \int_{-\infty}^{+\infty}D^h f(\tilde x -\Delta,v)\left(f(\tilde{x},v)-1\right)d\Delta dv\right)d \tilde{x}\]
				{Dropping the tildas, we can rewrite the inner term using the definition of $f$ as:
				\begin{align*}
					&\int_0^h \int_{-\infty}^{+\infty}D^h f(x-\Delta,v)(f(x,v)-1)d\Delta dv\\
					&=\int_0^h \int_{-\infty}^{+\infty} (f(x+\Delta,v)-f(x-\Delta,v))(f(x,v)-1)dvd\Delta\\
					&=-\int_0^h \int_{u(x)}^{u(x) \vee u(x+\Delta)\vee u(x-\Delta)} f(x+\Delta,v)-f(x-\Delta,v)dv d\Delta\\
					&= -\int_0^h \left(D^{\Delta} u\right)_+ -\left(D^{-\Delta} u\right)_+ d\Delta.
				\end{align*}
				We thus find 
				\[\overline{Q}_1=- \int_0^L\int_0^h \eta \left( \left(D^{\Delta} u\right)_+ -\left(D^{-\Delta} u\right)_+\right)  d\Delta dx.\]
				Using that 
				\[\int_0^L \eta \left(D^\Delta u\right)_+ dx=\int_0^L \eta(x-\Delta) \left(D^{-\Delta} u\right)_- dx,\]
				and similarly
				\[\int_0^L \eta \left(D^{-\Delta} u\right)_+ dx=\int_0^L \eta(x-\Delta) \left(D^{\Delta} u\right)_- dx,\]
we obtain
\begin{align*}
 \overline{Q}_1= & \frac{1}{2}\int_0^h\int_0^L  D^{\Delta} \eta \left(D^{\Delta} u\right)_+ -D^{-\Delta} \eta \left(D^{-\Delta} u\right)_+  +\eta^{-\Delta} D^{-\Delta} u -\eta^{\Delta} D^{\Delta} udx d\Delta\\
 =& \frac{1}{2} \int_0^h\int_0^L  D^{\Delta} \eta \left(D^{\Delta} u\right)_+ -D^{-\Delta} \eta \left(D^{-\Delta} u\right)_+ +  \eta (u(x-\Delta)-u(x+\Delta)) dx d\Delta\\
 =& \frac{1}{2}\int_0^h\int_0^L  D^{\Delta} \eta \left|D^{\Delta} u\right| +  \eta \left(u(x-\Delta)-u(x+\Delta)\right) dx d\Delta,
\end{align*}
which is \eqref{estimQ1}. If we now integrate \eqref{Q1} over time, we find
\begin{equation}\label{Q1plusQ2}
 Q_1+Q_2=\int_0^T\int_0^L\int_0^h D^\Delta \eta  |D^\Delta u|  d\Delta dx dt.
\end{equation}

				{\it Step 6.} We finally argue that 
				\begin{equation}\label{Q2}
				Q_3=-\left[\frac{1}{2} \int_0^h \int_0^L |D^\Delta u|D^\Delta u dx d\Delta \right]_0^T.
				\end{equation}
				Let us recall that by definition of $A$:					
				\begin{equation}\nonumber
					Q_3=\left[\int_{\R\times \R}\mathbb{1}_{\R_+}(v-w) \int_0^L \int_0^x D^h f(t,x,v)D^h f(t,y,w)dydxdvdw \right]_{t=0}^{t=T}.
				\end{equation}
				As above:
					\[
						Q_3=\left[- \int_0^L \int_0^{h}\int_{-\infty}^{+\infty}\int_{-\infty}^{v}   f(t,x,v)(f(t,x+\Delta,w)-f(t,x-\Delta,w)dwdvd\Delta dx\right]_{t=0}^{t=T}
						\]
and,
\[
 \int_{-\infty}^{v}  (f(t,x+\Delta,w)-f(t,x-\Delta,w)dw=(u(t,x+\Delta)-v\wedge 0)-(u(t,x-\Delta)-v \wedge 0). 
\]
Therefore
\begin{align*}
\int_{-\infty}^{+\infty}\int_{-\infty}^{v}   f(t,x,v)(f(t,x+\Delta,w)&-f(t,x-\Delta,w)dwdv\\
=&\int_{-\infty}^{u(t,x)}\int_{-\infty}^{v} (f(t,x+\Delta,w)-f(t,x-\Delta,w)dwdv\\
=&\int_{-\infty}^{u(t,x)}(u(t,x+\Delta)-v\wedge0)-(u(t,x-\Delta)-v\wedge 0) dv\\
=&\int_{(u(t,x+\Delta) \wedge u(t,x))}^{u(t,x)} u(t,x+\Delta)-v dv \\
&\qquad \quad- \int_{(u(t,x-\Delta)\wedge u(t,x))}^{u(t,x)} u(t,x-\Delta)-v dv
\end{align*}
 which, using that for $a,b \in \R$,
\[\int_{(a\wedge b)}^b (a-v) dv=-\frac{1}2 \left((a-b\wedge 0)\right)^2\]
						gives 
						\begin{align*}
						 Q_3&=\frac{1}{2}\left[ \int_0^L \int_0^{h}  (u(t,x+\Delta)-u(t,x)\wedge0)^2 -  (u(t,x-\Delta)-u(t,x)\wedge0)^2  d\Delta dx\right]_{t=0}^{t=T}\\
						 &=\frac{1}{2}\left[ \int_0^L \int_0^{h}  \left(D^{\Delta} u\right)_-^2 -  \left(D^{-\Delta}u\right)_-^2  d\Delta dx\right]_{t=0}^{t=T}\\
						 &=\frac{1}{2}\left[ \int_0^L \int_0^{h}  \left(D^{\Delta} u\right)_-^2 -  \left(D^{\Delta}u\right)_+^2  d\Delta dx\right]_{t=0}^{t=T}\\
						&=-\frac{1}{2}\left[ \int_0^L \int_0^{h}  |D^\Delta u| D^\Delta u  d\Delta dx\right]_{t=0}^{t=T}
						\end{align*}
which proves \eqref{Q2}. Combined with \eqref{Q1plusQ2}, this yields \eqref{KHMintequiv}.
}						\end{proof}

				We can now prove Proposition \ref{thm1}
				\begin{proof}[Proof of Proposition \ref{thm1}]
				By linearity, it is enough proving the estimate for $g=0$. Thanks to \eqref{KHMint} applied to $\eta= |\partial_x| \ksi$, we have
				 \begin{multline*} 
				 \int_0^T\int_0^L \left|D^h u\right|^3 dx dt \\= 6\int_0^T\int_0^L\int_0^h D^\Delta \eta   \left|D^\Delta u\right|  d\Delta dx dt-\left[3 \int_0^h \int_0^L  |D^\Delta u| D^\Delta u  dx d\Delta \right]_0^T.			  
				 \end{multline*}
 Thanks to Theorem \ref{equiv}, for $s<1$ and every function $v$, $\left\| |v| \right\|_{B^s_{p,r}}\le c \left\| v \right\|_{B^s_{p,r}}$ (notice that it also holds for $s=1$ if $p=r=2$). Using the triangle inequality and the invariance
 of the Besov norms with respect to translations, we obtain that for $s\in(0,1)$,  $\left\| |D^\Delta u| \right\|_{B^s_{p,r}}\le c \left\| u \right\|_{B^s_{p,r}}$. Applying \eqref{Besov1} we get
\begin{align*}\int_0^T\int_0^L\int_0^h D^\Delta \eta   \left|D^\Delta u\right|  d\Delta dx dt&=\int_0^T\int_0^L\int_0^h  \left(\left|\partial_x\right| D^\Delta \xi  \right) \left|D^\Delta u\right|  d\Delta dx dt \\
 &\le \frac{1}{\pi} \int_0^h  \left\| D^\Delta \ksi \right\|_{B^{1-s}_{p',r'}} \left\| \left|D^\Delta u\right| \right\|_{B^s_{p,r}}\\
 &\le c h \left\| \ksi \right\|_{B^{1-s}_{p',r'}} \left\| u \right\|_{B^s_{p,r}}.
\end{align*}
On the other hand, since
\[\left[\int_0^h \int_0^L  |D^\Delta u| D^\Delta u  dx d\Delta \right]_0^T\le c h \int_0^L u(0,x)^2+ u(T,x)^2 dx,\]
and since multiplying the equation \eqref{burger2} by $u$ and integrating gives,
					\begin{align*}
						\int_0^L \frac{1}{2}u(T,x)^2dx-\int_0^L \frac{1}{2}u(0,x)^2dx =&\int_0^T\int_0^L \partial_t\left(\frac{1}{2}u^2\right) dx dt\\
						=& \int_0^T\int_0^L u  \left|\partial_x \right|\ksi  dx dt \\
						\leq& c   \left\| \ksi \right\|_{B^{1-s}_{p',r'}}\left\| u \right\|_{B^s_{p,r}},
					\end{align*}
					we have
					\[\left[\int_0^h \int_0^L |D^\Delta u| D^\Delta u dx d\Delta \right]_0^T\le c h \left( \left\| \ksi \right\|_{B^{1-s}_{p',r'}}\left\| u \right\|_{B^s_{p,r}} + \|u(0,\cdot)\|_{L^2}^2\right).\] 
					Putting this together, we find 
					\[ \frac{1}{h}\int_0^T\int_0^L \left|D^h u\right|^3 dx dt\le c\left( \left\| \ksi \right\|_{B^{1-s}_{p',r'}}\left\| u \right\|_{B^s_{p,r}} + \|u(0,\cdot)\|_{L^2}^2\right)\]
					which concludes the proof.
				\end{proof}

				\begin{remark}
					Starting from \eqref{KHM}, one can obtain a larger family of estimates for the inhomogeneous Burgers equation. Unfortunately, the estimate
					\begin{align*}
					\left\| u \right\|_{B^{1/3}_{3,3}}^3  \leq&
								c \left(\left\|u\right\|_{B^{1/3}_{3,3}} \left\| g \right\|_{B^{2/3}_{3/2,3/2}} 
								+ \left\| u \right\|_{B^s_{p,r}} \left\| \ksi \right\|_{B^{1-s}_{p',r'}} +\left\| u(0,\cdot) \right\|_{L^2}^2  \right),
				\end{align*} 
				which would allow to avoid the logarithmic correction in \eqref{Main}, is borderline.
				\end{remark}

\section{Proof of Proposition \ref{thm3}}\label{sec:increaser}
				Let us remind the reader the statement we want to prove:
				\begin{proposition*}[Comparison between Besov norms of different index $r$] 
				There exists $c >0$ such that, for all $L>2$, for every $u$ solution of \eqref{Kuramoto} with average zero, the following estimate holds:
				\begin{align}\label{logpreuve}
					\left\|u \right\|_{\B^{1/3}_{3,3}} \leq c \ln^{1/3}(L)\left\|u\right\|_{\B^{1/3}_{3,\infty}}.
				\end{align}
			\end{proposition*}
						
			\begin{proof}
				The proof is elementary and resembles \cite[Prop. 4 II), Step 2 \&3]{Otto}. Let us first cut the term that we want to bound in three parts:
				\begin{align*}
					\int_0^{+\infty}  \frac{\left\|D^h  u \right\|^3_{L^3}}{h} \frac{dh}{h}=& \int_0^\ell \frac{\left\|D^h  u \right\|^3_{L^3}}{h} \frac{dh}{h}+ \int_\ell^L \frac{\left\|D^h  u \right\|_{L^3}^3}{h} \frac{dh}{h} + \int_L^{+\infty} \frac{\left\|D^h  u \right\|_{L^3}^3}{h} \frac{dh}{h}
					\\=&A(\ell)+B(\ell)+C.
				\end{align*}
				The large scale term  $C$ is in fact controlled by $A(\ell)$ and $B(\ell)$ by periodicity. Indeed:
				\begin{align*}
					C=& \int_L^{+\infty} \frac{\left\|D^h  u \right\|_{L^3}^3}{h} \frac{dh}{h}\\
					=&\sum_{n=1}^{+\infty} \int_{nL}^{(n+1)L}\frac{\left\|D^h  u \right\|^3_{L^3}}{h} \frac{dh}{h}\\
					\leq& \sum_{n=1}^{+\infty}\frac{1}{n^2} \int_0^L \frac{\left\|D^h  u \right\|^3_{L^3}}{h} \frac{dh}{h}\\
					\leq & c \left(A(\ell)+B(\ell)\right).
				\end{align*}
				The intermediate scale term $B(\ell)$ is directly handled with thanks to $\left\|u\right\|_{B^{1/3}_{3,\infty}}$:
				\begin{align*}
					B(\ell)=& \int_\ell^L \frac{\left\|D^h  u \right\|^3_{L^3}}{h} \frac{dh}{h}\\
					\leq & \sup_{h \in \R_+^*} \frac{\left\|D^h u\right\|^3_{L^3}}{h}  \int_\ell^L \frac{dh}{h}\\
					\leq & \ln(L/\ell) \left\|u\right\|_{B^{1/3}_{3,\infty}}^3.
				\end{align*}
				The small scale term $A(\ell)$ is the most difficult to bound. We shall prove that:
				\begin{equation}\label{pretendA}
					A(\ell) \leq c \ell  \left(\left(\int_0^L u(0,x)^2 dx \right)^{3/2} + L^{3/2} \left\| u\right\|^3_{B^{1/3}_{3,\infty}} \right).
				\end{equation}
				Before proceeding with the proof of \eqref{pretendA}, let us show that it is sufficient to conclude. Indeed, fix now $\ell=L^{-3/2}$. Then, we obtain:
				\begin{align*}
					\int_0^{+\infty}  \frac{\left\|D^h  u \right\|^3_{L^3}}{h} \frac{dh}{h} \leq c \ln(L) \left\|u\right\|_{B^{1/3}_{3,\infty}}^3 + c \left( \frac{1}{L}\int_0^L u(0,x)^2dx\right)^{3/2}.
				\end{align*}
				Taking now the time-space average yields  \eqref{logpreuve}. It remains to prove \eqref{pretendA}. \\
				
				 First, we prove that:
							\begin{equation}\label{changeRstep1}
				\int_0^T\int_0^L \left|D^h  u \right|^3dxdt
				\leq c h^2 \left( \left(\int_0^Lu(0,x)^2 dx\right)^{3/2}+\int_0^T\left(\int_0^L u^2dx \right)^{3/2}dt \right).
			\end{equation}
			We start by noting that:
			\begin{align*}
				\int_0^L \left|D^h  u(t,x)\right|^3dx
				 \leq& 2 \left(\sup_{x \in [0,L]}\left|u(t,x)\right|\right) \int_0^L \left| \int_x^{x+h} \partial_x u(t,y)dy \right|^2dx\\
				\leq& 2\left(\sup_{x \in [0,L]}\left|u(t,x)\right|\right) h^2 \int_0^L \left(\partial_x u \right)^2dx.
			\end{align*}
			For the convenience of the reader, we recall the argument for
			\begin{align}\label{majoresup}
				\sup_{x \in [0,L]}\left|u\right|
				\leq c\left(\int_0^L u^2 dx\right)^{1/4} \left(\int_0^L \left(\partial_x u \right)^2 dx\right)^{1/4}.
			\end{align}
			In fact, starting from:
			\begin{align*}
				u^2(t,x)= u^2(t,y)+ \int_y^x \frac{1}{2} u(t,z)\partial_xu(t,z)dz&& \forall y \in [0,L],
			\end{align*}
			and using that $u$ has zero average and thus vanishes somewhere, Cauchy-Schwarz inequality  immediately gives \eqref{majoresup}.
			Therefore,
			\begin{align*}
				\int_0^L \left|D^h  u\right|^3dx \leq& c h^2 \left(\int_0^L u^2dx \right)^{1/4} \left(\int_0^L \left(\partial_x u \right)^2dx \right)^{5/4},
			\end{align*}
			and using the following Sobolev inequality (which can easily be proved by using Fourier methods):
			\begin{align*}
				\int_0^L \left(\partial_xu\right)^2dx \leq \left(\int_0^L u^2dx \right)^{1/2}\left(\int_0^L \left(\partial_x^2u\right)^2dx \right)^{1/2},
			\end{align*}
			we get:
			\begin{align}\label{Sobo}
				\int_0^L \left|D^h  u\right|^3dx \leq& c h^2 \left(\int_0^L u^2dx \right)^{7/8} \left(\int_0^L \left(\partial_x^2 u \right)^2dx \right)^{5/8}.
			\end{align}
			Now, we have to work a bit to extract information from the energy identity. Multiplying \eqref{Kuramoto} by $u$ and integrating over $x$, we obtain:
			\begin{align*}
				\frac{d}{dt} \int_0^L u^2dx= &2 \int_0^L \left(\partial_xu\right)^2dx-2\int_0^L \left(\partial_x^2 u \right)^2dx\\
				\leq& 2 \left(\int_0^L u^2dx \right)^{1/2}\left(\int_0^L \left( \partial_x^2 u \right)^2dx \right)^{1/2}-2\int_0^L \left( \partial_x^2 u\right)^2dx \\
				\leq&  \int_0^L u^2dx-\int_0^L \left( \partial_x^2 u\right)^2dx.
			\end{align*}
			In a first step, we get from this differential inequality that there exists $c>0$ such that:
			\begin{align}\label{ENERGSTIM1}
				\int_0^L u^2(t+s,x)dx \leq c \int_0^L u^2(t,x)dx && \forall s \in [0,1].
			\end{align}
			In a second step, we get:
			\[
				\int_0^L u^2(t+1,x)dx - \int_0^L u^2(t,x)dx\leq  \int_t^{t+1} \int_0^L  u^2dxds-\int_t^{t+1} \int_0^L \left( \partial_x^2 u\right)^2dxds.
			\]
			The combination of both implies:
			\begin{equation}\label{ENERGSTIM2}
				\int_{t}^{t+1} \int_0^L \left(\partial_x^2u\right)^2 dxds + \sup_{s\in[t,t+1]} \int_0^L u^2(s,x) dx\leq c \int_0^L u^2(t,x)dx.
			\end{equation}
		Hence, together with \eqref{Sobo} in the form of
		\[\int_t^{t+1} \int_0^L  \left|D^h  u \right|^3 dx ds\le C h^2 \left( \sup_{s\in [t,t+1] }\int_0^L u^2(s,x) dx\right)^{7/8}\left(\int_t^{t+1}\int_0^L (\partial_x^2 u)^2 dx ds\right)^{5/8},\]
			we deduce
			\[\int_t^{t+1} \int_0^L  \left|D^h  u \right|^3 dx ds\le C h^2 \left(  \int_0^L u^2(t,x) dx\right)^{3/2}.\]
			Using this inequality for $t=0$ and in its integrated form between zero and $T$, we obtain \eqref{changeRstep1}.\\
		
		 We now claim that
			\begin{equation}\label{changeRstep2}
				\int_0^T\left(\int_0^L u^2 dx\right)^{3/2}dt \leq L^{3/2} \left\| u \right\|_{B^{1/3}_{3,\infty}}^3.
			\end{equation}
			Indeed, using H\"older's inequality, the fact that  $u$ has zero average and Jensen's inequality we obtain:
			\begin{align*}
				\int_0^T\left(\int_0^L u^2(t,x) dx\right)^{3/2} dt
				\leq& L^{1/2} \int_0^T\int_0^L \left|u(t,x)-\frac{1}{L}\int_0^L u(t,h)dh\right|^3dxdt\\
				\leq& L^{-1/2} \int_0^T \int_0^L \int_0^L \left|u(t,x)-u(t,x+h)\right|^3dh dx dt\\
				\leq& L^{1/2} \int_0^L \int_0^T \int_0^L \frac{\left|u(t,x+h)-u(t,x)\right|^3}{h}dxdtdh\\
				\leq& L^{3/2} \sup_{h \in \R_+} \int_0^T \int_0^L\frac{\left|D^h u(t,x)\right|^3}{h}dx dt.
			\end{align*}
		
			\subparagraph{Conclusion}
				Putting together \eqref{changeRstep1} and \eqref{changeRstep2}, we get as desired
				\begin{align*}
					A(\ell)=\int_0^\ell\int_0^T\int_0^L \left|D^h  u \right|^3dxdt \frac{dh}{h^2} \leq& c \int_0^\ell dh  \left( \left(\int_0^L u(0,x)^2 dx \right)^{3/2} + L^{3/2}\left\| u \right\|_{B^{1/3}_{3,\infty}}^3 \right)\\
					\leq& c \ell   \left( \left(\int_0^L u(0,x)^2dx \right)^{3/2} + L^{3/2}\left\| u \right\|_{B^{1/3}_{3,\infty}}^3 \right).
				\end{align*}
				
			\end{proof}
			\begin{remark}
			 As in \cite[Prop. 4 II), Step 2 \&3]{Otto}, we could have used an $L^\infty$ (in time) bound on $\int_0^L u^2 dx$ 
			 (proven for instance in \cite[Prop. 4]{Giacomelli}) to get \eqref{pretendA} directly. However, since we have a relatively simple and self-contained argument for it, we preferred to include it.
			\end{remark}
			\appendix
			
\section{Besov spaces}

		\subsection{Definition of time-space Besov spaces}
		We recall here some basics of the theory of  Besov spaces. We refer to \cite[Chapter 2, p. 51-121]{Bahouri}, for the construction of a dyadic Littlewood-Paley decomposition, and most of the proofs.
			\begin{definition}[Dyadic Littlewood-Paley decomposition]
				Let $\left(\phi_k\right)_{k\in \Z}$be a family of Schwartz functions such that their Fourier transforms $\left( \mathcal{F}\phi_k  \right)_{k \in \Z}$ satisfy:
				\begin{align*}
					&\mathcal{F}(\phi_0)(\ksi) = 0 &&\forall \left| \ksi \right| \notin \left]2^{-1}, 2\right[,  \\
					&\mathcal{F}(\phi_k)(\ksi)  = \mathcal{F}(\phi_0)\left(2^{-k}\ksi\right)& & \forall k \in \Z, \forall \ksi \in \R,	\\
					&\sum_{k \in \Z} \mathcal{F}(\phi_k)(\ksi)=1& & \forall \ksi \in \R. 
				\end{align*}
			\end{definition}
			Then, for a $L$-periodic function $u$, we define its Littlewood-Paley decomposition as:
			\begin{equation*}
				u_k(t,\cdot)= \phi_k * u(t,\cdot)
			\end{equation*}
			where $*$ denotes the periodic convolution.
			This allows us to define time-space Besov space $B^s_{p,r}$ for $s\in [0,+\infty]$, $p \in [1,\infty]$, $r\in [1,\infty]$ by the set of functions such that:
			\begin{align*}
				&\left\| u \right\|_{B^s_{p,r}}= \left( \sum_{k\in \Z} 2^{rsk} \left\| u_k \right\|^r_{L^p} \right)^{1/r} < \infty && \text{if } r<\infty,\\
				&\left\|u\right\|_{B^{s}_{p,r}}=\sup_{k\in\Z} 2^{sk} \left\|u\right\|_{L^p} < \infty && \text{if } r=\infty.
			\end{align*}
			We are actually interested in a rescaled homogeneous Besov norm, defined by:
			\begin{equation*}
				\left\| u \right\|_{\B^s_{p,r}}=\limsup_{T\rightarrow +\infty} \frac{1}{(LT)^{1/p}}\left\|u \right\|_{B^s_{p,r}}.
			\end{equation*}
			The time-space Besov norm can be replaced by an equivalent one, as stated in the following theorem (see \cite[Th. 2.36]{Bahouri}):
			\begin{theorem}\label{equiv} Let $s \in ]0,1[$ and $(p,r) \in [1,+\infty]^2$. Then there exists $c>0$ such that:
			\begin{equation*}
				c^{-1} \left\| \frac{\left\|D^h u\right\|_{L^P}}{h^s} \right\|_{L^r\left(\R_+,\frac{dh}{|h|}\right)} \leq \left\| u \right\|_{B^s_{p,r}} \leq c \left\| \frac{\left\|D^h u \right\|_{L^P}}{h^s} \right\|_{L^r\left(\R_+,\frac{dh}{|h|}\right)}.
			\end{equation*}
			where for $r<\infty$:
			\begin{equation*}
				\left\| \frac{\left\|D^h u\right\|_{L^P}}{h^s} \right\|_{L^r\left(\R_+,\frac{dh}{|h|}\right)}=\  \left(\int_0^{\infty} {\left(\int_0^T\int_0^L \left(\frac{\left|D^h  u \right|}{h^s}\right)^pdxdt \right)^{r/p}}\frac{dh}{h}\right)^{1/r},
			\end{equation*}
			and for $r=\infty$:
			\begin{equation*}
				\left\| \frac{\left\|D^h u\right\|_{L^P}}{h^s} \right\|_{L^{\infty}\left(\R_+,\frac{dh}{|h|}\right)}=\  \sup_{h>0} {\left(\int_0^T\int_0^L \left(\frac{\left|D^h  u \right|}{h^s}\right)^pdxdt \right)^{1/p}}.
			\end{equation*}
			\end{theorem}
			Besov spaces are particularly well adapted for interpolation as seen from the following theorem:
			\begin{theorem}[Interpolation between Besov spaces]\label{theointerpol}
				Let $(s,p,r),(s_1,p_1,r_1),(s_2,p_2,r_2) \in \R_+ \times [1,+\infty]^2$, and $u \in B^{s_1}_{p_1,r_1} \cap B^{s_2}_{p_2,r_2}$. If
				\begin{equation*}
				\left\{
					\begin{array}{r c l}
					s &=&\theta s_1 + (1-\theta) s_2, \\
					\frac{1}{p} &=& \frac{\theta}{p_1} + \frac{1-\theta}{p_2}, \\
					\frac{1}{r} &=& \frac{\theta }{r_1} + \frac{1-\theta}{r_2},
					\end{array}
				\right.
				\end{equation*}
				with $\theta \in ]0,1[$, then $u \in B^s_{p,r}$ and:
				\begin{equation}\label{Interpolation}
					\left\| u \right\|_{B^{s}_{p,r}} \leq \left\| u \right\|_{B^{s_1}_{p_1,r_1}}^{\theta} \left\| u \right\|_{B^{s_2}_{p_2,r_2}}^{1-\theta}.
				\end{equation}
			\end{theorem}				
			Theorem \ref{theointerpol} simply follows from the definition of the Besov norms and an application of H\"older's inequality. In some lemmas that we will enunciate later, we will use the partial derivative $\left|\partial_x \right|$ which is slightly different from the classical $\partial_x$. It is defined via Fourier series:
			\begin{equation}\nonumber
				\left|\partial_x \right|: \sum_{n\in\Z} a_n e^{\frac{2i\pi n x}{L}} \mapsto \sum_{n\in\Z}  a_n \frac{2\pi}{L} \left| n \right| e^{\frac{2i\pi n x}{L}}.
			\end{equation}
			The following theorems underlines a link between Besov spaces and the operator $\left|\partial_x \right|$:
			\begin{lemma}\label{prec}
				For all $\phi, g \in \mathcal{C}_L^{1}(\R)$:
				\begin{equation*}
					\int_0^L \phi \left| \partial_x\right|gdx =\frac{1}{\pi} \int_0^{+\infty} \int_0^L D^h  \phi D^h  g dx \frac{dh}{h^2}.
				\end{equation*}
			\end{lemma}
			\begin{proof}
				Let us expand $\phi$ and $g$ in Fourier series as
				\begin{align*}
					\phi(x)=\sum_{n\in \Z} \phi_n e^{\frac{2i\pi}{L}nx}  &&& \text{and}&&& g(x)=\sum_{n\in\Z} g_n e^{\frac{2i\pi}{L}nx}.
				\end{align*}
				Therefore, one can explicitly compute on the one hand:
				\begin{align*}
					\int_0^L \phi \left|\partial_x\right|g dx = \sum_{n\in\Z}2\pi \left|n\right| \phi_n  g_{-n},
				\end{align*}
				and on the other hand:
				\begin{align*}
					\int_0^{+\infty} \int_0^L D^h  \phi D^h  g dx \frac{dh}{h^2}=&
					\int_0^{+\infty} \sum_{n\in\Z} L\left( \left(e^{\frac{2i\pi}{L}h n}-1 \right) \phi_n \left(e^{-\frac{2i\pi}{L}h n}-1 \right)g_{-n} \right) \frac{dh}{h^2}\\
					=& \int_0^{+\infty} \sum_{n\in\Z}4L \sin^2\left(\frac{\pi h n}{L} \right) \phi_n  g_{-n} \frac{dh}{h^2}\\
					=& 4\sum_{n\in\Z} \phi_n g_{-n} \pi \left|n\right| \int_0^{+\infty} \sin^2\left(y\right) \frac{dy}{y^2}\\
					=& 2\pi^2 \sum_{n\in\Z} \phi_ng_{-n}|n|,
				\end{align*}
				which implies the result.
			\end{proof}
			We derive from this identity the following Besov estimate (see also \cite[Step. 3 p 39]{Otto}):
			\begin{proposition} \label{phidxg}
			Let $\phi, g \in \mathcal{C}^{1}_{L}(\R)$. Then, for all $p,p',r,' \in[1,+\infty]$, with $\frac{1}{p}+\frac{1}{p'}=1$ and $\frac{1}{r}+\frac{1}{r'}=1$, for all $s \in ]0,1[$ the following estimate holds:
				\begin{equation}\label{Besov1}
					\int_0^T \int_0^L \phi\left|\partial_x\right|g dx dt\leq \frac{1}{\pi} \left\| \phi \right\|_{B^s_{p,r}} \left\| g \right\|_{B^{1-s}_{p',r'}}.
				\end{equation}
			\end{proposition}
			\begin{proof}
				Using  Lemma \ref{prec} we get:
				\begin{align*}
					\int_0^T\int_0^L \phi \left|\partial_x\right| g dx dt &= \frac{1}{\pi} \int_0^{+\infty} \int_0^T \int_0^L D^h \phi D^h gdxdt\frac{dh}{h^2}.
				\end{align*}
				Then,  Hölder's inequality leads us to the result:
				\begin{align*}
					\int_0^{+\infty} \int_0^T \int_0^L D^h \phi D^h gdxdt\frac{dh}{h^2} \leq& \int_0^{+\infty} \left\|D^h \phi\right\|_{L^p} \left\|D^h  g \right\|_{L^{p'}} \frac{dh}{h^2}\\
					\leq& \left\| \frac{\left\|D^h \phi\right\|_{L^p}}{h^s} \right\|_{L^{r}\left(\R_+,\frac{dh}{h}\right)}
					 \left\| \frac{\left\|D^h g\right\|_{L^{p'}}}{h^{1-s}} \right\|_{L^{r'}\left(\R_+,\frac{dh}{h}\right)}.
				\end{align*}
			\end{proof}
	
			We finally state a useful lemma relating Besov norms of derivatives.
			\begin{lemma} \label{Lemme2.5}
				Let $s>0$, $p, r \in [1,\infty]$,  $m\in \mathbb{N}$ and $u \in B^{s+m-1}_{p,r}$. Suppose $h=\left| \partial_x \right|^{-1}\partial_x^mu$. Then there exists a positive constant $c$ depending only on $(s,p,r)$ such that the following estimate holds:
				\begin{equation}\label{transforme g ksi}
					\left\| h \right\|_{B^{s}_{p,r}} \leq c \left\| u \right\|_{B^{s+m-1}_{p,r}}.
				\end{equation}
			\end{lemma}
			
			\begin{proof}
				The proof is analogous to  \cite[Step 1 p. 17]{Otto}. By definition, we have:
				\begin{equation*}
					h_k = \phi_k * h.
				\end{equation*}
				Therefore, using the properties of convolution and of the quasi-orthogonality of the dyadic partition of unity we get:
				\begin{align*}
				h_k&=\phi_k*\sum_{k'\in[k-1;k+1]} \phi_{k'}*h\\
				&=\left|\partial_x \right|^{-1}\partial_x^m \phi_k * \sum_{k' \in [k-1,k+1]}u_{k'}.
				\end{align*}
				Then, using Young's inequality, we obtain:
				\begin{align*}
					\left\|h_k\right\|_{L^p} \leq& \left( \int_{\R} \left| \left| \partial_x\right|^{-1}\partial_x^m \phi_k\right|dx \right)  \sum_{k'\in[k-1,k+1]}\left\| u_{k'} \right\|_{L^p}\\
					\leq&2^{k(m-1)} \left( \int_{\R} \left| \left| \partial_x\right|^{-1}\partial_x^m \phi_0 \right|dx \right)\sum_{k'\in[k-1,k+1]}\left\| u_{k'} \right\|_{L^p}\\
					\leq& c 2^{k(m-1)}\sum_{k'\in[k-1,k+1]}\left\| u_{k'} \right\|_{L^p}.
				\end{align*}
				Hence:
				\begin{equation*}
					\sum_{k\in\Z} 2^{krs} \left\|h_k \right\|_{L^p}^r \leq c \sum_{k\in\Z} 2^{kr(m-1+s)}\left\|u_k \right\|_{L^p}^r,
				\end{equation*}
				which implies the aimed inequality.
			\end{proof}

		\section*{Acknowledgment}
		FO acknowledges many discussions with Dorian Goldman on simplifications of the proof of \cite{Otto}. MG was partially funded by a Von Humboldt fellowship. Most of this work was done while MG and MJ were hosted by the MPI-MIS Leipzig whose kind hospitality is acknowledged.
		\bibliographystyle{plain}
		\bibliography{Bib}
										
\end{document}